\newcommand{\si}[1]{#1}\newcommand{\jo}[1]{}

\si{
	\documentclass{article}
	\headheight=21.06892pt
	\addtolength{\textheight}{5cm}
	\addtolength{\topmargin}{-2.5cm}
	\setlength{\oddsidemargin}{-.4cm}
	\setlength{\textwidth}{17cm}
	\textheight=640pt
	\usepackage{doi}
}

\jo{
	\documentclass[moor]{informs1}              

	\usepackage{algorithm}
	\usepackage{algorithmic}
	\usepackage{tikz}
	\usetikzlibrary{shapes, shapes.geometric, arrows, positioning}
	\tikzstyle{nodo} = [rectangle, rounded corners, minimum width=1cm, minimum height=0.5cm,text centered, draw=black, fill=green!30]
	\tikzstyle{arrow} = [thick,->,>=stealth]
	
	\usepackage{natbib}
	 \NatBibNumeric
	 \def\newblock{\ }%
	 \bibpunct[, ]{[}{]}{,}{n}{}{,}%
	
	\usepackage[colorlinks=true,breaklinks=true,bookmarks=true,urlcolor=blue,
	     citecolor=blue,linkcolor=blue,bookmarksopen=false,draft=false]{hyperref}
	
	\def\EMAIL#1{\href{mailto:#1}{#1}}
	
	\TheoremsNumberedThrough     
	
	\EquationsNumberedThrough    
	
}

\usepackage{color}
\usepackage{graphicx}
\usepackage{multirow}
\usepackage{multicol}
\usepackage{mathtools}
\usepackage{tikz} 
\usetikzlibrary{calc,3d}
\usepackage{contour}
\usepackage{amsmath}
\usepackage[all,knot,arc,import,poly]{xy}
\usepackage{amsopn}
\usepackage{setspace}
\usepackage{todonotes}

\si{
	\usepackage{amsthm}
	\newtheorem{theorem}{Theorem}[section]
	\newtheorem{proposition}{Proposition}[section]
	\newtheorem{corollary}{Corollary}[section]
	\newtheorem{lemma}{Lemma}[section]
	\newtheorem{definition}{Definition}[section]
	\newtheorem{example}{Example}[section]
	\newtheorem{remark}{Remark}[section]
}

\newcommand{\inte}{\textnormal{int}}
\DeclareMathOperator{\bd}{bd}
\DeclareMathOperator{\bdp}{bd{^+}}

\newcommand{\R}{\mathbb{R}}

\newcommand{\N}{\mathbb{N}}
\newcommand{\linear}{\mathcal{L}_{\Omega}(\bar{x})}
\newcommand{\tangent}{\mathcal{T}_{\Omega}(\bar{x})}
\newcommand{\xb}{\bar{x}}

\renewcommand{\S}{\mathbb{S}}

\newcommand{\Ker}{\textnormal{Ker}}
\renewcommand{\Im}{\textnormal{Im}}

\renewcommand{\dim}{\textnormal{dim}}
\newcommand{\lin}{\textnormal{lin}}

\newcommand{\I}{\mathbb{I}}
\newcommand{\rank}{\textnormal{rank}}                
\newcommand{\seq}[1]{\{#1\}_{k\in \N}}
\newcommand{\E}{\mathbb{E}}
\newcommand{\K}{\mathcal{K}}
\newcommand{\cl}{\textnormal{cl}}
\newcommand{\cone}{\textnormal{cone}}
\newcommand{\A}{\mathcal{A}}
\newcommand{\V}{\mathcal{V}}
\newcommand{\C}{\mathcal{C}}
\newcommand{\G}{\mathcal{G}}
\newcommand{\viz}{\mathcal{N}}
\newcommand{\spn}{\textnormal{span}}
\newcommand{\pol}{^{\circ}}

\newcommand{\faceq}{\trianglelefteq}

\newcommand{\nlpcr}{CRCQ}
\newcommand{\nlprcr}{RCRCQ}
\newcommand{\socpcr}{CRCQ}
\newcommand{\sdpcr}{CRCQ}

\newcommand{\adj}{^T}

\usepackage{tikz}

\usepackage{amssymb,amsmath}
\usepackage{hyperref}
\hypersetup{
	colorlinks=true,
	citecolor=blue,
	linkcolor=blue,
	filecolor=magenta,      
	urlcolor=cyan,
}

\usepackage{latexsym}
\usepackage{mathtools}

\begin{document}

\si{
	\title{First- and second-order optimality conditions for second-order cone and semidefinite programming under a constant rank condition}
	\author{R. Andreani\footnote{Department of Applied Mathematics, University of Campinas, Campinas-SP, Brazil. Email: andreani@ime.unicamp.br} \and G. Haeser\footnote{Department of Applied Mathematics, University of S\~ao Paulo, S\~ao Paulo-SP, Brazil. Email: \{ghaeser,leokoto,thiagops\}@ime.usp.br}\and L. M. Mito\footnotemark[2] \and H. Ram\'{\i}rez C.\footnote{Departamento de Ingenier\'{\i}a Matem\'atica and Centro de Modelamiento Matem\'atico (CNRS UMI 2807), Universidad de Chile, Santiago, Chile. Email: hramirez@dim.uchile.cl}\and T. P. Silveira\footnotemark[2]}
}

\jo{
	\RUNAUTHOR{Andreani et al.}

	\RUNTITLE{Optimality conditions for NSOCP and NSDP under a constant rank CQ}
	
	\TITLE{First- and second-order optimality conditions for second-order cone and semidefinite programming under a constant rank condition}
	
	\ARTICLEAUTHORS{%
	\AUTHOR{Roberto Andreani}
	\AFF{Department of Applied Mathematics, University of Campinas, Campinas, Brazil, 
	\EMAIL{andreani@ime.unicamp.br}}
	
	\AUTHOR{Gabriel Haeser, Leonardo M. Mito, Thiago P. Silveira}
	\AFF{Department of Applied Mathematics, University of São Paulo, São Paulo, Brazil,
	\EMAIL{ghaeser@ime.usp.br}, 
	\EMAIL{leokoto@ime.usp.br},
	\EMAIL{thiagops@ime.usp.br}}
	
	\AUTHOR{H{\'e}ctor Ram{\'i}rez C.}
	\AFF{Departamento de Ingenier\'{\i}a Matem\'atica and Centro de Modelamiento Matem\'atico (CNRS UMI 2807), Universidad de Chile, Santiago, Chile, 
	\EMAIL{hramirez@dim.uchile.cl}}	
	} 
	
	\ABSTRACT{%
	The well known constant rank constraint qualification [Math. Program. Study 21:110--126, 1984] introduced by Janin for nonlinear programming has been recently extended to a conic context by exploiting the eigenvector structure of the problem. In this paper we propose a more general and geometric approach for defining a new extension of this condition to the conic context. The main advantage of our approach is that we are able to recast the strong second-order properties of the constant rank condition in a conic context. In particular, we obtain a second-order necessary optimality condition that is stronger than the classical one obtained under Robinson's constraint qualification, in the sense that it holds for every Lagrange multiplier, even though our condition is independent of Robinson's condition.
	}%
	

	\MSCCLASS{Primary: 90C46; Secondary: 90C30, 90C22.}
	
}

\maketitle

%
%
	
\si{
	\abstract{The well known constant rank constraint qualification [Math. Program. Study 21:110--126, 1984] introduced by Janin for nonlinear programming has been recently extended to a conic context by exploiting the eigenvector structure of the problem. In this paper we propose a more general and geometric approach for defining a new extension of this condition to the conic context. The main advantage of our approach is that we are able to recast the strong second-order properties of the constant rank condition in a conic context. In particular, we obtain a second-order necessary optimality condition that is stronger than the classical one obtained under Robinson's constraint qualification, in the sense that it holds for every Lagrange multiplier, even though our condition is independent of Robinson's condition.}

\vspace{0.3cm}

}

%
%
\section{Introduction}

In the classical \textit{nonlinear programming} (NLP) context, the so-called \textit{constant rank constraint qualification} (\nlpcr{})~\cite{crcq} was first presented as a tool for stability analysis, which stood out for being independent of the usual \textit{Mangasarian-Fromovitz constraint qualification} (MFCQ) and strictly weaker than the \textit{linear independence constraint qualification} (LICQ). For instance, it has been applied with this purpose in NLP~\cite{tilt,crcq,rcrcq,minch,tiltcrmf}, \textit{mathematical programs with equilibrium constraints} (MPEC)~\cite{sensitivitympec}, \textit{generalized equations}~\cite{outrata1}, and \textit{bilevel optimization}~\cite{mehlitzminch,xuye}. Also, it is the origin of several other constant rank-type conditions, such as the \textit{constant positive linear dependence}~\cite{rcpld,ams05,cpld} and the \textit{constant rank of the subspace component}~\cite{cpg}, which have been successfully applied in the convergence analysis of iterative algorithms. To name a few algorithms whose convergence theory relies on CRCQ and its variants, we point out: an augmented Lagrangian method~\cite{abms,ams2}, a regularized interior point method~\cite{qiu}, sequential quadratic programming methods for NLP~\cite{sqpliu,cpld,sqpwright} and MPEC~\cite{sqpjiang}, and some relaxation schemes for MPEC~\cite{relaxmpec,ulbrich}. In fact, a particularly interesting aspect of CRCQ that makes it suitable for supporting practical algorithms is the fact it can be roughly interpreted as a relaxation of LICQ that is able to separate the core information of the problem, ignoring redundant constraints. Moreover, all linear programming problems satisfy CRCQ, in contrast with LICQ and MFCQ.

Besides convergence of algorithms and stability analysis, CRCQ was used in several contexts, such as NLP~\cite{aes2010,ams2,minchleschov}, MPEC~\cite{secondordermpec}, vector optimization~\cite{vectoropt}, and continuous-time NLP~\cite{continuousnlp}, for studying necessary second-order optimality conditions. One of the main goals of this paper is to bring such results to more general conic programming contexts, namely \textit{nonlinear second-order cone programming} (NSOCP) and \textit{nonlinear semidefinite programming} (NSDP). As far as we know, the best second-order results for these problems have appeared in a well-known paper by Bonnans, Cominetti, and Shapiro~\cite{bonn-comi-shap}. Basically, they derived no-gap second-order optimality conditions for problems over \textit{second-order regular cones}~\cite[Definition 3]{bonn-comi-shap}, such as NSDP and NSOCP, under the well-known \textit{Robinson's CQ} (see~\eqref{conic:rcq} on page~\pageref{nondegeneracycq}, or~\cite{robinson}), which is the natural extension of MFCQ to conic programming. In particular, their second-order necessary condition states that every local solution that satisfies Robinson's CQ must also satisfy the following: for every critical direction, there exists a Lagrange multiplier (possibly depending on this direction), such that a certain quadratic form is nonnegative with respect to such direction and multiplier. However, the second-order condition that is obtained under \nlpcr{} in NLP replaces ``there exists a Lagrange multiplier'' with ``for every Lagrange multiplier'', which is stronger than the one of~\cite{bonn-comi-shap}. Although this stronger condition can be obtained from~\cite{bonn-comi-shap} after assuming that the Lagrange multiplier is unique, which is ensured by stronger constraint qualifications such as the \textit{nondegeneracy} condition (see  \eqref{nondegeneracycq} on page \pageref{nondegeneracycq}), this assumption is often regarded as too stringent. To the best of our knowledge, no second-order result concerning every Lagrange multiplier, without assuming its uniqueness, has been presented so far in the literature of nonlinear conic programming. Moreover, no extension of \nlpcr{} has been proposed for nonlinear conic programming until very recently.

In 2019, Zhang and Zhang~\cite{ZZ} proposed an extension of \nlpcr{} and its relaxed version~\cite{rcrcq} for NSOCP, but it was later discovered that their results were incorrect~\cite{andreani2020erratum}. This event has motivated us to investigate other possible extensions of \nlpcr{} to conic problems, and their properties. The first step in this direction was made in~\cite{andreani2020naive}, for NSOCP and NSDP problems with multiple constraints. The idea of~\cite{andreani2020naive} is to rewrite some of the conic constraints as locally equivalent NLP constraints, whenever possible, and then jointly applying nondegeneracy and the NLP version of \nlpcr{} to the resulting problem. Later, in~\cite{seqcrcq}, based on the ideas from~\cite{weak-sparse-cq}, we improved this strategy by exploiting the eigenvector structure of the semidefinite cone to deal with the conic constraints that could not be rewritten as NLP constraints. This approach was also extended to NSOCP problems in~\cite{seqcrcq-socp}. In simple terms, the condition of~\cite{seqcrcq,seqcrcq-socp} demands the rank of some families of functions to remain constant along every sequence converging to the point of interest -- roughly speaking, a constant rank ``by paths'' -- therefore, this extension is highly specialized to deal with sequences generated by iterative algorithms, but since this rank may vary between paths, it is likely unsuitable for other purposes. Indeed, the focus of~\cite{seqcrcq,seqcrcq-socp} was the global convergence of a large class of algorithms to first-order stationary points, and no second-order results were provided in it. Nevertheless, it is reasonable to expect that~\nlpcr{} may have multiple independent and correct extensions, each one of them generalizing at least one important aspect of it, but perhaps not all of them.

A common feature of all previous attempts of extending CRCQ to a conic context is an approach based on re-characterizing the conic program and the nondegeneracy condition, trying to make them as similar to NLP and LICQ as possible, so the extension of CRCQ would come out straightforwardly. This is somehow understandable because, even in NLP, the CRCQ condition has never received a geometrical interpretation before. In this paper, we present a new geometrical characterization of \nlpcr{} for NLP in terms of the faces of the nonnegative orthant, which suggests a natural extension of it to NSOCP and NSDP. A point that we should stress is that contrary to our previously mentioned works, the definition of CRCQ that we present here is very simple. We prove that this extension is a constraint qualification strictly weaker than nondegeneracy and independent of Robinson's CQ, as it should be, and we also compare it with the condition of~\cite{seqcrcq,seqcrcq-socp}. Then, as an application, we show that every local solution of the problem satisfies the strong second order optimality condition, provided our extension of \nlpcr{} holds. Moreover, just as it happens in NLP, our result does not demand \textit{a priori} any specific condition over the Lagrange multiplier set, besides nonemptiness.

The structure of this paper is as follows: Section~\ref{section_conic} consists of a nonlinear conic programming review emphasizing some aspects of the theory that are not commonly discussed in the literature; in Section~\ref{sec:nlp}, we analyze \nlpcr{} for NLP and we show how it can be interpreted in terms of the faces of the nonnegative orthant. In Sections~\ref{section_socp} and~\ref{sec:sdp}, we propose extensions of \nlpcr{} for NSOCP and NSDP, respectively, and we prove some of its properties. Finally, in Section~\ref{sec:conc}, we conclude this paper with a short discussion and some ideas of prospective work.

We end this section by introducing some of our basic notation:
throughout this paper, $\mathbb{E}$ will denote a finite-dimensional linear space equipped with the inner product $\langle \cdot, \cdot \rangle$; and for a given set $S\subseteq \E$, we will denote the \emph{polar} of $S$ by
\[
	S\pol:=\{z\in \E\mid \langle z, y\rangle\leq 0, \ \forall y\in S\}\]
and the \textit{orthogonal complement} of $S$ will be denoted by $S^\perp$. The notations $\cl(S)$, $\inte(S)$, $\bd(S)$, and $\bd^+(S)$ stand for the topological closure, interior, boundary, and boundary excluding the origin of $S$ in $\E$, respectively. The \textit{smallest cone} that contains $S$ will be denoted by $\cone(S)$, and the \textit{smallest linear space} that contains $S$ will be denoted by $\spn(S)$. Moreover, we denote by $Dg(x)$ the derivative of a twice continuously differentiable function $g:\R^n\to\mathbb{E}$ at a given point $x\in \R^n$, and by $Dg(x)\adj$ the adjoint of $Dg(x)$, which by definition satisfies $\langle Dg(x)d, z \rangle = \langle d, Dg(x)\adj z \rangle$ for all $d\in \R^{n}$ and $z \in \mathbb{E}$. Similarly, $D^2g(x)$ denotes the second-order derivative of $g$ at $x$, and the action of $D^2 g(x)$ over $d_1,d_2\in\R^n$ will be denoted by $D^2g(x)[d_1,d_2]$.

%
%
\section{Common framework: nonlinear conic programming}\label{section_conic}

In this section, we will review some classical results of convex analysis, and first- and second-order optimality conditions and constraint qualifications for NSOCP and NSDP. These problems are the cornerstones of two independent research fields, but they can also be seen as particular cases of a \textit{nonlinear conic programming} (NCP) problem, given by
\begin{equation}
\begin{array}{ll}
\nonumber\mbox{\textnormal{Minimize }} 	& f(x),\\ 
\label{conic_problem}\mbox{s.t. }	& g(x)\in \mathcal{K},
\end{array}
\tag{NCP}
\end{equation}
where $f:\R^{n} \to \R$ and $g:\R^{n} \to \mathbb{E}$ are twice continuously differentiable, and $\mathcal{K}\subseteq \mathbb{E}$ is a closed convex pointed cone that is assumed to be nonempty. We will use~\eqref{conic_problem} as a framework to discuss the common traits of NSOCP and NSDP simultaneously, before moving to specific traits. Throughout the whole paper, we will denote the feasible set of~\eqref{conic_problem} by $\Omega := \{x \in \R^{n} \mid g(x) \in \mathcal{K} \}$. 

Let us begin with two key ideas that underlie all the results of this paper: \textit{reducibility} and \textit{faces}. Recall from~\cite[Definition 3.135]{bonnans-shapiro} that for any given linear spaces $\E$ and $\mathbb{F}$, a cone $\K\subseteq \E$ is said to be \textit{reducible} (more precisely, $C^2$-\emph{reducible}) at a point $y\in \K$, to a closed convex pointed cone $\mathcal{C}\subseteq \mathbb{F}$, if there exists a neighborhood $\viz$ of $y$ and a twice continuously differentiable reduction function $\Xi\colon \viz\to \mathbb{F}$ (possibly depending on $y$) such that $\Xi(y)=0$, $D\Xi(y)$ is surjective, and
\[
	\K\cap \viz=\{z\in \viz\mid \Xi(z)\in \mathcal{C}\}.
\]
In general, reductions are meant to be used as a simplification tool that allows one to interpret any point of $\K$ as a vertex of some other cone $\C$, and then extend the results obtained at $\C$ to $\K$ in a smooth way. In this work, we are also interested in the geometrical properties of the reduced cone $\mathcal{C}$ as well; in particular, in its faces.

To make a brief revision, we recall that $F$ is a \emph{face} of $\mathcal{C}$ if every open line segment that contains a point of $F$ also has its extrema in $F$; that is, if for every $y\in F$ and every $z,w\in \mathcal{C}$ such that $y=\alpha z + (1-\alpha)w$ for some $\alpha\in (0,1)$, we have that $z,w\in F$. Further, when there exists some $\eta\in \mathcal{C}\pol$ such that 
\[F=\mathcal{C}\cap\{\eta\}^\perp,\]
that is, when $F$ is the intersection between $\mathcal{C}$ and one of its supporting hyperplanes, we say that $F$ is an \emph{exposed face} of $\C$. Some cones, like the nonnegative orthant, the semidefinite cone, and the second-order cone, are \emph{facially exposed}, meaning all of their faces are exposed. We use the notation $F\faceq \mathcal{C}$ to say that $F$ is a face of $\mathcal{C}$. 

Now, to contextualize our results, we will revisit the classical theory of NCP in the next section, with a special emphasis in the work of Guignard~\cite{guignard}, and Bonnans, Cominetti, and Shapiro~\cite{bonn-comi-shap}. In particular, we stress some aspects of the NCP theory that are often disregarded in the literature.

\subsection{Review of first-order optimality conditions}
For any set $S\subseteq \E$ and any $z\in S$, recall the (Bouligand) \textit{tangent cone} to $S$ at $z$, defined as
\begin{displaymath}
\mathcal{T}_{S}(z) := \left\{y \in \mathbb{E}\mid \exists \seq{t_k} \rightarrow 0^{+}, \exists \seq{y^k} \rightarrow y \textnormal{ such that } z + t_{k}y^{k} \in S \textnormal{ for all } k\in \N\right\}.
\end{displaymath}

Our review of first-order constraint qualifications for~\eqref{conic_problem} revolves around two particular cones: the tangent cone $\mathcal{T}_{\Omega}(\xb)$ to $\Omega$ at a feasible point $\xb\in \Omega$, and the \textit{linearized tangent cone}
\begin{displaymath}
\linear := \left\{d \in \R^{n} \mid Dg(\xb)d \in \mathcal{T_{K}}(g(\xb))\right\},
\end{displaymath}
where $\mathcal{T}_{\K}(g(\xb))$ is the tangent cone to $\K$ at $g(\xb)$. The importance of these cones for our analyses lies on the necessary optimality conditions for~\eqref{conic_problem} associated with them.  Namely, given any local minimizer $\xb\in \Omega$ of \eqref{conic_problem}, it is easy to see that $\langle\nabla f(\xb), d\rangle \geq 0$ for all $d \in \tangent$; that is,
\begin{equation}\label{def:geometric_condition}
	-\nabla f(\xb) \in \tangent^{\circ}.
\end{equation}
This is one of the simplest necessary optimality conditions, sometimes called the \textit{first-order geometric necessary condition} for the optimality of $\xb$. However, it may be difficult to use~\eqref{def:geometric_condition} when $\Omega$ does not admit an explicit characterization since $\tangent\pol$ may not be easily computable in this case. The polar of $\linear$, on the other hand, admits a practical description, as it is shown in the following lemma, extracted from the proof of~\cite[Theorem 2]{guignard} by Guignard: 
\begin{lemma}\label{closedness}
Let $\xb \in \Omega$. Then, $\linear\pol=\cl (H(\xb))$, where 
\begin{equation}\label{h(x)}
H(\xb) := Dg(\xb)\adj \mathcal{N_{K}}(g(\xb))=\left\{ Dg(\xb)\adj z\mid z\in \mathcal{N}_{\K}(g(\xb))\right\},
\end{equation}
and $\mathcal{N}_{\K}(g(\xb)):=\mathcal{T_{K}}(g(\xb))\pol$ is the \emph{normal cone} to $\K$ at $g(\xb)$.
\end{lemma}
\jo{\proof{Proof.}}\si{\begin{proof}}
By the bipolar theorem (see e.g. \cite[Proposition 2.40]{bonnans-shapiro}), it suffices to prove that $\linear=H(\xb)\pol$. Take any direction $d \in \linear$ and let $z \in \mathcal{T_{K}}(g(\xb))^{\circ}$. By definition, $Dg(\xb)d \in \mathcal{T_{K}}(g(\xb))$ and then
\[
	0\geq\langle Dg(\xb)d, z\rangle = \langle d, Dg(\xb)\adj z \rangle.
\] Thus, since $z$ is arbitrary, we obtain that $d \in H(\xb)^{\circ}$; and since $d$ is also arbitrary, it follows that $\linear \subseteq \left(H(\xb)\right)^{\circ}$. Conversely, assume that there exists a vector $v \in H(\xb)^{\circ}$ such that $v \notin \linear$, that is, $Dg(\xb)v \notin \mathcal{T_{K}}(g(\xb))$. By the strong separation theorem (see e.g. \cite[Theorem 2.14]{bonnans-shapiro}), there exists a vector $y$ such that $\langle y, Dg(\xb)v \rangle > 0$ and $\langle y, z\rangle < 0$, for all $z \in \mathcal{T_{K}}(g(\xb))$, that is, $y \in \mathcal{N_{K}}(g(\xb))$. Therefore, $Dg(\xb)\adj y\in H(\xb)$, which is a contradiction with $\langle Dg(\xb)\adj y, v \rangle > 0$, since $v\in H(\xb)^{\circ}$.
\jo{\hfill\Halmos \\}

\si{\end{proof}}

Recall that since $\K$ is a closed convex cone, we have
\[
	\mathcal{N}_{\K}(g(\xb))=\left\{ z\in \K\pol\mid \langle g(\xb),z\rangle=0 \right\}.
\]
Then, combining the first-order geometric necessary condition and Lemma~\ref{closedness} yields the following theorem, also by Guignard:
\begin{theorem}[Theorem 2 of~\cite{guignard}]\label{kkt:guignard}
Let $\xb\in \Omega$ be a local minimizer of \eqref{conic_problem}. If $\tangent\pol=\linear\pol$ and $H(\xb)$ is closed, then there exists some $\bar\lambda\in \K\pol$ such that
\begin{equation}\label{conic_kkt}
\nabla f(\xb) + Dg(\xb)\adj \bar\lambda = 0
 \ \textnormal{ and } \
\langle g(\xb), \bar\lambda \rangle = 0.
\end{equation}
\end{theorem}

Theorem~\ref{kkt:guignard} can be seen as the ``dual form'' of the first-order geometric condition~\eqref{def:geometric_condition}, and any vector $\bar\lambda\in \K\pol$ that satisfies the \textit{Karush-Kuhn-Tucker conditions}~\eqref{conic_kkt} is called a \textit{Lagrange multiplier} associated with $\xb$. Moreover, the collection of all Lagrange multipliers associated with $\xb$ will be denoted by $\Lambda(\xb)$, and when $\Lambda(\xb)\neq \emptyset$ we say that $\xb$ is a \textit{KKT point} of~\eqref{conic_problem}. 

The hypothesis of Theorem~\ref{kkt:guignard}, 
\begin{equation}\label{guignardcq}
	\tangent\pol=\linear\pol \ \textnormal{ and } \ H(\xb) \textnormal{ is closed,}
\end{equation}
is known in the literature as \emph{Guignard's CQ}, and it is the weakest assumption that makes the KKT conditions necessary for the local optimality of $\xb$, in the sense of: if $\Lambda(\xb)\neq \emptyset$ for every continuously differentiable function $f$ that has a local minimizer constrained to $\Omega$ at $\xb$, then Guignard's CQ must also hold at $\xb$~\cite[Corollary 3.4]{gould1975optimality}. B\"orgens et al.~\cite[Definition 5.11]{borgens2020new} defined Guignard's CQ for optimization problems in Banach spaces as a single equality
\[
	\tangent^{\circ} = H(\xb),
\]
which is equivalent to~\eqref{guignardcq} due to Lemma~\ref{closedness}. In NLP, Guignard's CQ is usually stated in the form $\tangent\pol=\linear\pol$, since the closedness of $H(\xb)$ follows from the polyhedricity of $\R^m_+$. However, as it can be seen in the following example, the equality $\tangent\pol=\linear\pol$ on its own may not ensure that $\Lambda(\xb)\neq \emptyset$ when $H(\xb)$ is not closed.

\begin{example}\label{socp_example1}
Consider the following problem, presented in~\cite[Subsection 2.1]{andersen2002notes}:
\begin{equation}\label{socp_example1_problem}
\begin{array}{ll}
\nonumber\mbox{\textnormal{Minimize }} 	& f(x):=-x_2,\\ 
\mbox{s.t. }	& g(x):=(x_1, x_1, x_2) \in K_{3},
\end{array}
\end{equation}
where $K_{3}$ is the three-dimensional second-order cone, given by \[K_{3} = \left\{(x_{1},x_{2},x_{3}) \in \R^{3}\mid x_{1} \geq \sqrt{x_{2}^{2} + x_{3}^{2}}\right\}.\] Note that its feasible set is given by $\Omega = \{x\in \R^{2}\mid x_1 \geq 0 \textnormal{ and } x_2=0 \}$, and that the point $\xb = (0,0) \in \R^{2}$ is a local minimizer of it. Any Lagrange multiplier $\lambda:=(\lambda_1,\lambda_2,\lambda_3)\in K_3^\circ$ associated with $\xb$ must satisfy
\begin{equation}\label{socp_example1_kkt}
\left(\begin{array}{c} 0 \\ -1 \end{array} \right) +\lambda_{1}\left(\begin{array}{c} 1 \\ 0 \end{array}\right) +\lambda_{2} \left(\begin{array}{c} 1 \\ 0 \end{array}\right) + \lambda_{3}\left(\begin{array}{c} 0 \\ 1 \end{array}\right) = \left(\begin{array}{c} 0 \\ 0 \end{array}\right),
\end{equation}
which implies that $\lambda_{3} = 1$ and $\lambda_{1} = -\lambda_{2}$. But since $\lambda \in K_{3}^\circ=-K_{3}$, then this vector must also satisfy $-\lambda_{1} \geq \sqrt{\lambda_{1}^{2} + 1}$, which does not have a  solution with $\lambda_{3} = 1$ and $\lambda_{1} = -\lambda_{2}$. Therefore, $\xb$ does not satisfy the KKT conditions. However, note that $\tangent=\Omega=\linear$ and consequently, $\tangent\pol=\linear\pol$. Additionally, note that
\[
	H(\xb)=\{(y_1+y_2,y_3)\in \R^2\mid (y_1,y_2,y_3)\in K_3\pol\}
\]
is not closed, because the sequence $\left\{\left(-\frac{1}{k},-1\right)\right\}_{k\in \N}$ is contained in $H(\xb)$ since
	$\left(-\frac{1}{k}-k, k, -1 \right)\in K_3\pol, \ \forall k\in \N,$
but its limit point $(0,-1)$ does not belong to $H(\xb)$.
\jo{\hfill\Halmos \\}

\end{example}
The condition 
\begin{equation}\label{abadiecq}
	\tangent=\linear  \ \textnormal{ and } \ H(\xb) \textnormal{ is closed,}
\end{equation}
which implies Guignard's CQ, is known as \emph{Abadie's CQ} (see also~B\"orgens et al.~\cite[Definition 5.5]{borgens2020new}), and Example~\ref{socp_example1} tells us that the closedness of $H(\xb)$ cannot be omitted in this case, either. The reason why we emphasize this point is that, as far as we know, it appears that Abadie's CQ and Guignard's CQ are rarely seen in the literature of finite-dimensional conic programming problems other than NLP, and the closedness of $H(\xb)$ is rarely regarded in the study of constraint qualifications. In contrast, $H(\xb)$ plays an important role in our results.

In finite-dimensional conic contexts, the focus is usually on constraint qualifications that already imply $H(\xb)$ is closed without requiring it explicitly, such as Robinson's CQ, that holds at a given point $\xb\in \Omega$ when
\begin{equation}\label{conic:rcq}
0\in \inte(\Im (Dg(\xb))-\K+g(\xb)).
\end{equation}
In particular, if $\K$ has nonempty interior, then Robinson's CQ holds at $\xb$ if, and only if, there exists some $d\in \R^n$ such that
\[
	g(\xb)+Dg(\xb)d\in \inte(\K).
\]
Robinson's CQ is stronger than Abadie's CQ, and it implies that $\Lambda(\xb)$, besides being closed and convex, is also nonempty and bounded~\cite[Theorem 3.9]{bonnans-shapiro} when $\xb$ is a local minimizer of \eqref{conic_problem}. 
If $\K$ is reducible at the point $g(\xb)$ to a cone $\mathcal{C}$ by the reduction function $\Xi$, then the constraint $g(x)\in\K$ is locally equivalent to the reduced constraint $\mathcal{G}(x)\in \mathcal{C}$, where
$	\mathcal{G}:=\Xi\circ g.$
In this case, Robinson's CQ holds at $\xb$ for the original constraint if, and only if, it holds for the reduced constraint at the same point.

Another well-known constraint qualification in the context of conic programming is the \textit{nondegeneracy} condition, which holds at $\xb$ when
\begin{equation}\label{nondegeneracycq}
	\Im(Dg(\xb))+\lin(\mathcal{T}_{\K}(g(\xb)))=\E,
\end{equation}
where $\lin(\mathcal{T}_{\K}(g(\xb)))=\mathcal{T}_{\K}(g(\xb))\cap -\mathcal{T}_{\K}(g(\xb))$ denotes the largest linear space contained in $\mathcal{T}_{\K}(g(\xb))$; that is, its \emph{lineality space}. This CQ has first appeared in Shapiro and Fan's article~\cite{shapiro1995eigenvalue} for NSDP, by the name \textit{transversality}, and then it was generalized to NCP by Shapiro, in~\cite{shapiro-uniqueness}. Nondegeneracy is strictly stronger than Robinson's CQ and it is known that if $\xb$ is a local minimizer of~\eqref{conic_problem} that satisfies nondegeneracy, then $\Lambda(\xb)$ is a singleton (see, for instance,~\cite[Proposition 4.75]{bonnans-shapiro}). Moreover, if $\K$ is reducible, nondegeneracy is equivalent to the surjectivity of $D\mathcal{G}(\xb)$, as it can be easily deduced from the equality $\lin(\mathcal{T}_{\K}(g(\xb)))=\Ker(D \Xi(g(\xb)))$; see \cite[Section 4.6.1]{bonnans-shapiro}.

Due to their implications over the Lagrange multiplier set, nondegeneracy and Robinson's CQ are currently the most important CQs in the study of second-order optimality conditions for~\eqref{conic_problem}, which will be reviewed in the next subsection.

\subsection{Second-order optimality conditions}

Before starting, recall that the (inner) second-order tangent set to a nonempty set $S\subseteq\E$, at a point $z\in S$, in a direction $y\in \mathcal{T}_S(z)$, is defined by
\begin{equation}
	\label{def:tan2}
	\begin{aligned}
		\mathcal{T}^2_{S}(z,y) & :=\left\{w\in \E \ \middle| \ z+ty+\frac{t^2}{2}w+o(t^2)\in S, \ \forall t> 0\right\},\\
	\end{aligned}
\end{equation}
which is closed for all such $z$, $y$, and $S$. In addition, if $S$ is convex, then $\mathcal{T}^2_{S}(z,y)$ is also convex~\cite[Page 163]{bonnans-shapiro}; and if $S$ is second-order regular, as it is the case of the semidefinite cone and the second-order cone, then $\mathcal{T}^2_S(z,y)$ is nonempty~\cite[Page 202]{bonnans-shapiro}. 

The role of second-order necessary optimality conditions is to provide additional information when first-order conditions are not meaningful enough; that is, along the directions in the cone
\[
	C(\xb):=\left\{ d\in \R^n\mid d\in \tangent, \ \langle\nabla f(\xb), d\rangle= 0\right\},
\]
which is often called the \textit{cone of critical directions}, or simply, the \textit{critical cone} of~\eqref{conic_problem} at $\xb$. Ben-Tal and Zowe~\cite{ben1982unified} presented a \textit{geometric second-order necessary optimality condition} for~\eqref{conic_problem}, stating that if $\xb$ is a local minimizer of the problem, then
\begin{equation}\label{geom2}
	\langle\nabla f(\xb), s\rangle + \langle \nabla^2 f(\xb)d, d\rangle\geq 0
\end{equation}
for every $d\in C(\xb)$ and every $s\in \mathcal{T}^2_{\Omega}(\xb,d)$.  Then, Kawasaki~\cite[Theorem 5.1]{kawasaki1988envelope} made the first advances to derive a ``dual form'' of~\eqref{geom2} under Robinson's CQ assuming that $\K$ is a closed convex cone with nonempty interior. This result was later generalized and refined by Cominetti~\cite[Theorem 4.2]{cominetti1990metric} to the case where $\K$ is assumed to be a closed convex set. An important improvement was made afterwards by Bonnans, Cominetti, and Shapiro~\cite{bonn-comi-shap}, who clarified several key points of the previous works, and obtained no-gap\footnote{The term ``zero gap'', or ``no gap'', is often used in NLP to refer to a second-order condition that does not require constraint qualifications to be necessary (using Fritz John/generalized Lagrange multipliers), and that becomes sufficient after replacing an inequality by a strict inequality. However, in this paper, we say that a condition has zero gap when it satisfies the latter, possibly subject to a constraint qualification, in the same way as~\cite{bonn-comi-shap}.} second-order conditions, in particular, for second-order regular cones~\cite[Section 4]{bonn-comi-shap}.
Let us recall Bonnans, Cominetti, and Shapiro's necessary condition in the context of second-order regular cones:

\begin{theorem}[Theorem 3.1 of~\cite{bonn-comi-shap}]\label{conic:bsoc}
Let $\xb\in \Omega$ be a local minimizer of~\eqref{conic_problem} that satisfies Robinson's CQ. Then, for every direction $d\in C(\xb)$, there exists some $\bar\lambda_d\in \Lambda(\xb)$, such that
\begin{equation}
d^T \nabla^2 f(\xb)d+\langle D^2g(\xb)[d,d], \bar\lambda_d \rangle - \sigma(d,\xb,\bar\lambda_d)\geq 0,
\end{equation}
where
\begin{equation}\label{conic:sigma}
	\sigma(d,\xb,\bar{\lambda}_d):=\sup\left\{ \langle w, \bar{\lambda}_d\rangle\mid w\in \mathcal{T}^2_{\K}(g(\xb),Dg(\xb)d) \right\}
\end{equation}
is the support function of $\mathcal{T}^2_{\K}(g(\xb),Dg(\xb)d)$ with respect to $\bar{\lambda}_d$. 
\end{theorem}

The term $\sigma(d,\xb,\bar{\lambda}_d)$ characterizes a possible curvature of the set $\K$ at $g(\xb)$ along $Dg(\xb)d$, and it is often called the ``sigma-term'' in the classical literature (for instance, in the book~\cite{bonnans-shapiro}). 
Since $\bar\lambda_d\in \Lambda(\xb)$ and $\K$ is convex, $\sigma(d,\xb,\bar{\lambda}_d)$ is always nonnegative; and if $\K$ is polyhedral, as in NLP, then the sigma-term is zero everywhere. See also the discussion on polyhedricity and extended polyhedricity in~\cite[Section 3.2.3]{bonnans-shapiro}. It is also worth mentioning that the second-order optimality condition of Theorem~\ref{conic:bsoc} can be derived without constraint qualifications, using Fritz John (generalized) multipliers~\cite[Theorem 3.50]{bonnans-shapiro}.

Although the condition of Theorem~\ref{conic:bsoc} is generally considered very natural and useful in the conic programming context and in NLP, a stronger condition where the Lagrange multiplier $\bar{\lambda}$ does not depend on $d$ has several potential uses, in view of the NLP literature. This motivates the following definition:

\begin{definition}\label{conic:ssoc}
Let $\xb\in \Omega$ be a KKT point and let $\bar\lambda\in \Lambda(\xb)$ be given. We say that the pair $(\xb,\bar{\lambda})$ satisfies the \emph{second-order condition} (SOC) when
\begin{equation}\label{conic:eqsoc}
	d^T \nabla^2 f(\xb)d+\langle D^2g(\xb)[d,d],  \bar\lambda \rangle - \sigma(d,\xb,\bar\lambda)\geq 0,
\end{equation}
for every $d\in C(\xb)$.
\end{definition}

In NLP, the existence of some $\bar{\lambda}\in \Lambda(\xb)$ such that SOC holds for the pair $(\xb,\bar\lambda)$ is known as the \textit{semi-strong second-order necessary optimality condition}~\cite{bonnans-semistrong}. Moreover, when SOC holds for every $\bar{\lambda}\in \Lambda(\xb)$, then we obtain what is known as the~\textit{strong second-order necessary optimality condition}~\cite{aes2010}. However, while the condition of Theorem~\ref{conic:bsoc} is necessary for optimality under Robinson's CQ, this is not true, in general, for the strong and semi-strong conditions. In fact, there is a counterexample published by Baccari~\cite[Section 3]{Baccari2004} (see also Anitescu~\cite{anitescu} and Arutyunov~\cite{arutyunov1998second}), that shows that Robinson's CQ does not guarantee the existence of a $\bar{\lambda}\in \Lambda(\xb)$ such that the pair $(\xb,\bar{\lambda})$ satisfies SOC (see also the extended version of~\cite{ninoconjecture} for details). Under nondegeneracy, the set $\Lambda(\xb)$ is a singleton and, in this case, the semi-strong and the strong second-order conditions both coincide with the condition of Theorem~\ref{conic:bsoc}. 

As far as we know, there is no result concerning the semi-strong and strong second-order conditions without assuming uniqueness of Lagrange multipliers in the literature of conic programming, except for NLP. In NLP, this has been addressed by means of constant rank-type constraint qualifications, which is also the path we will follow in this paper.

%
%

\section{Revisiting constant rank CQs in NLP}\label{sec:nlp}

In this section we will revisit some constant rank-type conditions for NLP from a geometrical point of view, in order to extend it to a more general conic context later on. Consider the standard NLP problem
\begin{equation}\label{NLP}
\begin{array}{lll}
\nonumber\mbox{\textnormal{Minimize }} 	& f(x) , &\\ 
\mbox{s.t. }	& g_j(x) \geq 0, & j=1,\ldots,m,\\
& g_j(x) = 0, & j=m+1,\ldots,m+p,
\end{array}
\tag{NLP}
\end{equation}
which is a particular case of \eqref{conic_problem} with $\E=\R^{m+p}$, $\K=\R^m_+\times \{0\}^p$, and $g(x):= (g_1(x),\ldots,g_{m+p}(x))$. As usual in NLP, given a feasible point $\xb$ of \eqref{NLP}, we will denote the set of active inequality constraints at $\xb$ as $\A(\xb) := \{j \in \{1,\ldots,m\}\mid g_{j}(\xb) = 0\}$.

Now, let us recall Janin's constant rank constraint qualification as it was first presented in~\cite{crcq}.

\begin{definition}[\nlpcr{} \cite{crcq}]\label{nlp_crcq_definition}
Let $\xb$ be a feasible point of (\ref{NLP}). We say that the \textnormal{constant rank constraint qualification for NLP} (\nlpcr{}) holds at $\xb$ if there exists a neighborhood $\V$ of $\xb$ such that, for every subset $J \subseteq \A(\xb)\cup \{m+1,\ldots,m+p\}$, the rank of the family $\{\nabla g_{j}(x)\}_{j \in J}$ remains constant for all $x \in \V$.
\end{definition}

To prove that \nlpcr{} is a constraint qualification, Janin proved that it implies $\linear \subseteq \tangent$, which in turn implies Abadie's CQ in NLP. His proof is what motivates the requirement to consider every subset $J$ of $\A(\xb)\cup\{m+1,\ldots,m+p\}$ in Definition~\ref{nlp_crcq_definition}; indeed, after picking a direction 
\[
	d \in \linear = 
	\left\{d \in \R^{n} \left|
		\begin{array}{l}
			\nabla g_j(\xb)^{T}d \geq 0, \ j \in \A(\xb),\\
			\nabla g_j(\xb)^{T}d = 0, \ j \in \{m+1,\ldots,m+p\}
		\end{array}	
		\right.
	\right\},
\]
in order to prove that $d\in \tangent$, it is sufficient to have the constant rank assumption for the constraints that correspond to the indices $j\in \mathcal{A}(\xb)$ such that $\nabla g_{j}(\xb)^{T}d = 0$. Since those indices depend on $d$, and they are not determined \textit{a priori}, one considers all possibilities. However, as it was noted several years later by Minchenko and Stakhovski~\cite{rcrcq}, taking subsets of the equality constraints is completely superfluous, even for Janin's proof. The ``correct'' definition of \nlpcr{} was then presented in~\cite{rcrcq} as a relaxed version of \nlpcr{}.

\begin{definition}[\nlprcr{}~\cite{rcrcq}]
Let $\xb$ be a feasible point of (\ref{NLP}). We say that \textnormal{relaxed constant rank constraint qualification for NLP} (\nlprcr{}) holds at $\xb$ if there exists a neighborhood $\V$ of $\xb$ such that, for every subset $J \subseteq \A(\xb)$, the rank of the family $\{\nabla g_{j}(x)\}_{j \in J\cup \{m+1,\ldots,m+p\}}$ remains constant for all $x \in \V$.
\end{definition}

In order to bring these CQs to the conic setting, 
our approach in this manuscript consists first in 
generalizing two key ideas of NLP: the notion of ``active constraints'' and the notion of ``subsets of indices of active constraints''. The former can be interpreted in the general context as a consequence of reducibility. Indeed, for any given $\xb\in \Omega$, let $s:=|\mathcal{A}(\xb)|$ and note that $\R^m_+\times \{0\}^p$ is reducible at $g(\xb)$ to the cone 
\[
	\mathcal{C}:=\R^{s}_+\times \{0\}^p
\]
in a neighborhood $\viz$ of $g(\xb)$ by the mapping $\Xi\colon \viz\to \R^{s+p}$ such that \[
	\Xi(y):=(y_j)_{j\in \mathcal{A}(\xb)\cup\{m+1,\ldots,m+p\}}\]
for every $y\in \viz$, and in this case the reduced constraint function of~\eqref{NLP} at $\xb$ takes the form
\begin{equation}\label{nlp:reducedcons}
	\mathcal{G}(x):=\Xi(g(x))=(g_j(x))_{j\in \mathcal{A}(\xb)\cup\{m+1,\ldots,m+p\}}.
\end{equation}
Therefore, in NLP, reducing the problem is essentially the same as simply disregarding inactive constraints around the point $\xb$. The notion of ``subsets of indices of the active constraints'', on the other hand, can be interpreted in terms of faces.

It is easy to see that every face of $\R^s_+$ can be written in terms of a unique subset of the canonical vectors of $\R^s$, which we will denote by $c_1,\ldots,c_s$. That is, $F\faceq \R^s_+$ if, and only if, there exists some $J\subseteq \{1,\ldots,s\}$ such that
\begin{equation}\label{nlp:eqface}
	F=\R^s_+\bigcap_{j\in J} \{c_i\}^\perp,
\end{equation}
where $F$ and $J$ are clearly in a one-to-one correspondence.

\si{\newpage}
\tikzset{cart/.style={x={(1cm,0cm)}, y={(-1.25cm,-0.75cm)}, z={(0cm,1cm)}}}
\begin{figure}[!ht]\label{fig:faces}
\centering
\begin{tikzpicture}[scale=0.5,inner sep=0.2cm]
  \begin{scope}[cart]
   \draw[-latex, thick, color=black] (0,0,0) -- (4,0,0) node[anchor=west] {$c_1$};
   \draw[-latex, thick, color=black] (0,0,0) -- (1.5,-1.5,0) node[anchor=north] {$c_2$};
   \draw[-latex, thick, color=black] (0,0,0) -- (0,0,4) node[anchor=east] {$c_3$};
   \fill[color=red, opacity=0.3] (0,0,0) -- (4,0,0) -- (5,-1.5,0) -- (1.5,-1.5,0) -- cycle;
   \fill[color=green, opacity=0.3] (0,0,0) -- (1.5,-1.5,0) -- (1.5,-1.5,4) -- (0,0,4) -- cycle;
   \fill[color=blue, opacity=0.3] (0,0,0) -- (4,0,0) -- (4,0,4) -- (0,0,4) -- cycle;
  \end{scope}
   \draw[ball color=black] (0,0,0) circle (0.05);
\end{tikzpicture}
 \caption{Faces of $\R^3_+$}
\end{figure}

For example, in Figure 1, the vertex of $\R^3_+$ corresponds to $J=\{1,2,3\}$; the one-dimensional faces $\cone(c_1)$, $\cone(c_2)$, and $\cone(c_3)$ correspond to $J=\{2,3\}$, $J=\{1,3\}$, and $J=\{1,2\}$, respectively; the left, front, and bottom two-dimensional faces correspond to $J=\{1\}$, $J=\{2\}$, and $J=\{3\}$, respectively; and $\R^3_+$ itself corresponds to $J=\emptyset$. 

Thus, considering all subsets of active constraints at $\xb$ is the same as considering all faces of the reduced cone $\mathcal{C}=\R^{s}_+\times \{0\}^p$. This discussion suggests a natural characterization of \nlprcr{} in terms of the faces of the reduced cone, as follows:

\begin{proposition}\label{nlp:crcq-conic-improved}
Let $\xb$ be a feasible point of (\ref{NLP}). Then, \nlprcr{} holds at $\xb$ if, and only if, there exists a neighborhood $\V$ of $\xb$ such that, for each $F\faceq \R^{|\mathcal{A}(\xb)|}_+\times \{0\}^p$, the dimension of
\[
	D\mathcal{G}(x)^T [F^\perp]
\]
remains constant for every $x\in \V$, where $\mathcal{G}$ is as defined in~\eqref{nlp:reducedcons}.
\end{proposition}

\jo{\proof{Proof.}}\si{\begin{proof}}
Let $s:=|\mathcal{A}(\xb)|$ and, without loss of generality, let us assume that $\mathcal{A}(\xb)=\{1,\ldots,s\}$. Moreover, let $c_1,\ldots,c_{s+p}$ be the canonical basis of $\R^{s+p}$, and let $F\faceq \R^{s}_+\times \{0\}^p$. Note that $F=R\times\{0\}^p$, where $R\faceq\R^{s}_+$. Then, there exists some $J\subseteq\{1,\ldots,s\}$ such that
\begin{equation*}\label{nlp:eqface2}
	F=\left(\R^s_+\bigcap_{j\in J} \{c_i\}^\perp\right) \times \{0\}^p,
\end{equation*}
which implies
\[
	F^\perp=R^\perp\times\R^p=\textnormal{span}\left(\{c_j\mid j\in J\cup\{s+1,\ldots,s+p\}\}\right),
\]
so
\begin{equation}\label{nlp:rcrcqeq}
	D\mathcal{G}(x)^T [F^\perp]=\textnormal{span}(\{D\mathcal{G}(x)^T c_j\}_{j\in J\cup \{s+1,\ldots,s+p\}})=\textnormal{span}(\{\nabla g_j(x)\}_{j\in J\cup \{m+1,\ldots,m+p\}}).
\end{equation}
Consequently,
\[
	\dim(Dg(x)^T [F^\perp])=\textnormal{rank}(\{\nabla g_j(x)\}_{j\in J\cup \{m+1,\ldots,m+p\}}).
\]
The conclusion follows from the one-to-one correspondence between $F$ and $J$.
\jo{\hfill\Halmos \\}

\si{\end{proof}}
The equivalent form of \nlprcr{} presented in Proposition~\ref{nlp:crcq-conic-improved} allows us to visualize what it actually describes, geometrically. Indeed, recall that $\R^n=D\G(x)^{-1}(\spn(F))+(D\G(x)^{-1}(\spn(F)))^\perp$ and it is elementary to see that
\[
	(D\G(x)^{-1}(\spn(F)))^\perp=D\G(x)\adj[F^\perp].
\]
This implies the following relation:
\[
	\dim(D\G(x)^{-1}(\spn(F)))+\dim(D\G(x)\adj[F^\perp])=n.
\]
Thus, \nlprcr{} can be equivalently stated as the constant dimension of $D\G(x)^{-1}(\spn(F))$ for every $x\in \mathcal{V}$ at each $F\faceq \C=\R^{|\mathcal{A}(\xb)|}_+\times \{0\}^p$. The set $D\G(x)^{-1}(\spn(F))$, on the other hand, can be regarded as a ``linear approximation'' of $\G^{-1}(\C)$ around $\xb$, since $D\mathcal{G}(x)$ is the best linear approximation of $\mathcal{G}$ at $x\in \V$ and, similarly, the faces of $\mathcal{C}$ can also be seen as ``linear approximations'' of it at $\G(\xb)$. In fact, each face induces a potentially different linear approximation of $\G^{-1}(\C)$, which in turn coincides with $\Omega$ around $\xb$. So roughly speaking: \nlprcr{} holds at $\xb$ when the dimension of every linear approximation of the feasible set $\Omega$ at $\xb$ is invariant to small perturbations. In particular, defining $g_{J}(x):=(g_j(x))_{j\in J\cup \{m+1,\ldots,m+p\}}$ for every $J\subseteq \mathcal{A}(\xb)$, this characterization is equivalent to the constant dimension of $\Ker(Dg_J(x))$ for all $x$ in a neihborhood of $\xb$ at every $J\subseteq \mathcal{A}(\xb)$, which can also be trivially seen from the original definition of~\nlprcr{}.

Note that the characterization of \nlprcr{} from Proposition~\ref{nlp:crcq-conic-improved} and the discussion above do not appear to be limited to the context of NLP, contrary to its original definition. In the next two sections, we will prove that the same idea can be applied to NSOCP and NSDP, respectively, giving rise to new constraint qualifications.

\begin{remark}\label{nlp:remcrcq}
It is possible to obtain a characterization of \nlpcr{} in the same style of Proposition~\ref{nlp:crcq-conic-improved}. To do this, it suffices to reformulate the equality constraints $g_j(x)=0$ as a pair of inequality constraints $g_j(x)\geq 0$ and $-g_j(x)\geq 0$, for $j\in \{m+1,\ldots,m+p\}$. That is, consider $\K:=\R^m_+\times \R^p_+\times \R^p_+$ and $g(x):=(g_1(x),\ldots,g_{m+p}(x),-g_{m+1}(x),\ldots,-g_{m+p}(x))$ in Proposition~\ref{nlp:crcq-conic-improved}.
\end{remark}

In view of Remark~\ref{nlp:remcrcq}, we see that there are multiple ways of dealing with equality constraints in our approach, and they are not all equivalent. The suitability of each approach may depend on the application, but we highlight that our approach is able to deal with equality constraints regardless of how they are modelled. For simplicity, equality constraints are omitted in our exposition. See also Remarks~\ref{socp_remark_eq} and~\ref{sdp_remark_eq}. In the following two sections, we extend the ideas of this section to NSOCP and NSDP.

%
%
\section{Nonlinear second-order cone programming}\label{section_socp}

In this section, we consider the following problem:
\begin{equation}\label{socp_multifold}
\begin{array}{lll}
\nonumber\mbox{\textnormal{Minimize }} 	& f(x), &\\ 
\mbox{s.t. }	& g_{j}(x)\in K_{m_{j}}, & j=1,\ldots, q,
\end{array}
\tag{NSOCP}
\end{equation}
where $K_{m_{j}} := \{(z_{0}, \widehat{z}) \in \R \times \R^{m_{j}-1}\mid z_{0} \geq \|\widehat{z}\| \}$ when $m_{j} > 1$ and $K_{1} = \{x \in \R\mid x \geq 0\}$. 
Since $K_{m_j}$ is self-dual, we have that $z\in K_{m_j}\pol$ if, and only if, $-z\in K_{m_j}$, for any $j=1,\ldots,q$. Also, note that~\eqref{socp_multifold} can be seen as a particular case of~\eqref{conic_problem} with 
\[
	\K:=K_{m_1}\times\ldots\times K_{m_q} \ \textnormal{ and } \ g(x):=(g_1(x),\ldots,g_q(x)).
\]

Given a feasible point $\xb\in \Omega$, let us define the following index sets:
\[
	\begin{aligned}
	I_{\inte}(\xb) & := \{j\ \in \{1,\ldots,q\} \mid g_{j}(\xb) \in \inte (K_{m_j})\},\\
I_{B}(\xb) & := \{j \in \{1,\ldots,q\} \mid g_{j}(\xb) \in \bdp(K_{m_j})\},\\
I_{0}(\xb) & := \{j \in \{1,\dots,q\} \mid g_{j}(\xb)=0\},
	\end{aligned}
\]
which consist of the indices 
of the constraints that hit the interior, the boundary excluding zero, and the vertex of their respective cones. For simplicity, we will omit equality constraints; we should mention, nevertheless, that our results can be easily adapted to deal with equality constraints --- see Remark~\ref{socp_remark_eq} for details. As another measure to avoid cumbersome notation, we will assume that $I_B(\xb)=\{1,\ldots,|I_B(\xb)|\}$; this assumption will often be recalled throughout this section.

Following Bonnans and Ram{\'i}rez~\cite{BonRam}, for any given $\xb\in \Omega$, we see that $\K$ is reducible to 
\begin{equation}\label{socp:defC}
	\mathcal{C}:=\prod_{j\in I_0(\xb)} K_{m_j}\times \R^{|I_B(\xb)|}_+
\end{equation}
in a neighborhood $\viz_1\times\ldots\times\viz_q$ of $g(\xb)$ by the function $\Xi:=(\Xi_j)_{j\in I_0(\xb)\cup I_B(\xb)}$, where $\Xi_j\colon \viz_j\to \R^{m_j}$ is the identity function for every $j\in I_0(\xb)$, and $\Xi_j\colon \viz_j\to \R$ is given by
\begin{equation}\label{socp:reductionmap}
\Xi_{j}(y) := y_{0} - \|\widehat{y}\|
\end{equation}
for every $j\in I_B(\xb)$, and every $y\in \R^{m_j}$.  This leaves us with the reduced constraint
\[
	\mathcal{G}(x)\in \mathcal{C},
\]
where $\mathcal{G}(x):=\Xi(g(x))=(\mathcal{G}_j(x))_{j\in I_0(\xb)\cup I_B(\xb)}$,
\begin{equation}\label{socp:defG}
	\mathcal{G}_j(x):=\Xi_j(g_j(x))=\left\{
		\begin{array}{ll}
			g_j(x), & \textnormal{ if } \ j\in I_0(\xb),\\ 			\phi_j(x), & \textnormal{ if } \ j\in I_B(\xb),
		\end{array}			
	\right.
\end{equation}
and $\phi\colon \R^n\to \R^{|I_B(\xb)|}$ has its $j$-th component given by
\begin{equation}\label{socp:reductionmap2}
	\phi_{j}(\xb) :=[g_{j}(x)]_{0} - \|\widehat{g_{j}(\xb)}\|.
\end{equation}
Note that $g(x)\in \K$ if, and only if, $\mathcal{G}(x)\in \mathcal{C}$ for every $x$  sufficiently close to $\xb$.

By~\cite[Lemma 25]{BonRam}, we see that the linearized cone of the original constraints of~\eqref{socp_multifold} at a given $\xb\in \Omega$ can be computed as
\begin{displaymath}
\linear = \left\{ \begin{tabular}{c|cl}
\multirow{2}{*}{$d \in \mathbb{R}^{n}$} & $Dg_{j}(\xb)d \in K_{m_{j}},$    & $j \in I_{0}(\xb)$ \\ 
                                        & $ D\phi(\xb) d \in \R^{|I_B(\xb)|}_+$    & \\ 
\end{tabular}  \right\},
\end{displaymath}
and that it coincides with the linearized cone of the reduced constraint at $\xb$. Moreover, it follows from~\cite[Lemma 15]{AG} that for each $j=I_{\inte}(\xb)\cup I_B(\xb)$, we have $\langle \bar\lambda_{j}, g_{j}(\xb) \rangle = 0$, if, and only if,
\begin{equation}\label{socp:multformat}
\bar\lambda_{j} =
\left\{
\begin{array}{ll}
	0, & \textnormal{ if } j \in I_{\inte}(\xb),\\
	\frac{[\bar\lambda_j]_0}{[g_j(\xb)]_0} R_{m_{j}}g_{j}(\xb), & \textnormal{ if } j \in I_{B}(\xb),
\end{array}
\right.
\end{equation}
where $R_{m_j}$ is a matrix defined as
\begin{equation}
	R_{m_{j}}:=\begin{bmatrix}
	1 & 0\\
	0 & -\I_{m_j-1}
	\end{bmatrix},
\end{equation}
and $\I_{m_j-1}$ is the $(m_j-1)\times (m_j-1)$ identity matrix. Therefore, still following~\cite{BonRam}, the point $\xb$ satisfies the KKT conditions with respect to the constraint $g(x)\in \K$ if, and only if, there exist some vectors $\bar{\lambda}_{j} \in K_{m_{j}}\pol$, $j\in I_0(\xb)\cup I_B(\xb)$, such that:
\begin{equation}\label{socp_kkt_reduction}
\nabla f(\xb) + \sum_{j \in I_{0}(\xb)} Dg_{j}(\xb)^{T}\bar\lambda_{j} + \sum_{j \in I_{B}(\xb)} \frac{[\bar\lambda_j]_0}{[g_j(\xb)]_0} Dg_j(\xb)^T R_{m_j}g_j(\xb) = 0,
\end{equation}
which also coincides with the KKT conditions with respect to the reduced constraint $\mathcal{G}(x)\in \mathcal{C}$. In fact, note that for each $j\in I_B(\xb)$, the reduced Lagrange multiplier with respect to the reduced constraint $\phi_j(x)\geq 0$ is simply $[\bar\lambda_j]_0$.

With this in mind, we are ready to present our extension of \nlpcr{} (and \nlprcr{}) to NSOCP inspired by the characterization of Proposition~\ref{nlp:crcq-conic-improved}.

%
%
\subsection{A facial constant rank constraint qualification for NSOCP}

Recall that, for each $j=1,\ldots,q$, the cone $K_{m_j}$ is facially exposed, meaning every $F\faceq K_{m_j}$ can be written as the intersection of one of its supporting hyperplanes, say $\{\eta\}^\perp$ with $\eta\in K_{m_j}$. In fact, although $K_{m_j}$ has infinitely many faces when $m_j>2$, they are limited to only three types: 
\begin{itemize}
\item The vertex, $\{0\}$, which can be characterized by any $\eta\in \inte(K_{m_j})$;
\item The cone $K_{m_j}$ itself, which is characterized by $\eta=0$;
\item A ray at the boundary of $K_{m_j}$, starting at the vertex and passing through a point $z\in \bd^+(K_{m_j})$, which can be characterized by any $\eta\in \cone(R_{m_j}z)\setminus \{0\}$.
\end{itemize}

Moreover, every $F\faceq \mathcal{C}$ has the form
\[
	F=\left(\prod_{j\in I_0(\xb)} F_j\right)\times R,
\]
where $F_j\faceq K_{m_j}$ for every $j\in I_0(\xb)$, and $R\faceq \R^{|I_B(\xb)|}_+$. Then, for every $x\in \R^n$, sufficiently close to $\xb$, we have
\begin{displaymath}
	D\mathcal{G}(x)\adj[F^\perp]=\sum_{j\in I_0(\xb)}Dg_{j}(x)^T [F_j^{\perp}] + D\phi(x)^T[R^\perp],
\end{displaymath}
where $\phi(x):=(\phi_j(x))_{j\in I_B(\xb)}$. This motivates the following definition:

\begin{definition}\label{socp_facial_condition}
Let $\xb$ be a feasible point of~\eqref{socp_multifold}. We say that the \emph{facial constant rank property} holds at $\xb$ if there exists a neighborhood $\V$ of $\xb$ such that for each $F\faceq \mathcal{C}$, the dimension of $D\mathcal{G}(x)\adj[F^\perp]$ remains constant for all $x \in \V$, where $\G$ is given by~\eqref{socp:defG} and $\C$ is given by~\eqref{socp:defC}.
\end{definition}

Recall the discussion after Proposition~\ref{nlp:crcq-conic-improved} and note that Definition~\ref{socp_facial_condition} can be equivalently stated in terms of the constant dimension of $D\G(x)^{-1}(\spn(F))$ for all $x\in \mathcal{V}$ and every $F\faceq \C$. That is, the facial constant rank property holds at $\xb$ when the dimension of every linear approximation of the feasible set remains locally invariant around $\xb$. Although this characterization is somewhat more intuitive than Definition~\ref{socp_crcq_definition}, the latter is easier to use.

The facial constant rank condition is sufficient for the equality $\tangent = \linear$ to hold. 
To prove this, we employ the main result of Janin's paper~\cite{crcq}, but the version we use is a slightly different characterization found in~\cite[Proposition 3.1]{aes2010}. Despite the fact we work in a context more general than NLP, we use the same result that was used in NLP. 

\begin{proposition}(\cite[Proposition 3.1]{aes2010})\label{prop1} Let $\{\zeta_{i}(x)\}_{i \in \mathcal{I}}$ be a finite family of twice continuously differentiable functions $\zeta_i\colon \R^n\to \R$, $i\in \mathcal{I}$, such that the family of its gradients $\{\nabla \zeta_{i}(x)\}_{i \in \mathcal{I}}$ remains with constant rank in a neighborhood of $\xb$, and consider the linear subspace 
\[
	\mathcal{S} := \left\{y \in \R^{n}\mid \langle \nabla \zeta_{i}(\xb), y\rangle = 0,  \ i \in \mathcal{I}\right\}.
\]
Then, there exists some neighborhoods $V_{1}$ and $V_{2}$ of $\xb$, and a diffeomorphism $\psi: V_{1} \rightarrow V_{2}$, such that:

\begin{itemize}
\item[(i)] $\psi(\xb) = \xb$;
\item[(ii)] $D\psi(\xb)=\mathbb{I}_n$;
\item[(iii)] $\zeta_{i}(\psi^{-1}(\xb+y))=\zeta_{i}(\psi^{-1}(\xb))$ for every $y \in \mathcal{S}\cap (V_2-\xb)$ and every $i\in \mathcal{I}$.
\end{itemize}
Moreover, the degree of differentiability of $\psi$ is the same as of $\zeta_{i}$, for all $i\in \mathcal{I}$.
\end{proposition}

For the last part of the above proposition, about the degree of differentiability of $\psi$, we refer to Minchenko and Stakhovski~\cite[Page 328]{minch}. Now, we are able to prove the main result of this section:

\begin{theorem}\label{socp_facial_theorem}
Let $\xb$ be a feasible point of (\ref{socp_multifold}). If the facial constant rank property holds at $\xb$, then $\tangent = \linear$.
\end{theorem}

\jo{\proof{Proof.}}\si{\begin{proof}}
It suffices to show that $\linear\subseteq\tangent$. Let $d\in\linear$ and suppose that $\xb$ satisfies the facial constant rank property. Let 
\begin{equation}\label{def:minface-full}
	F:=\left(\prod_{j\in I_0(\xb)} F_j\right)\times R,
\end{equation}
where $F_j\faceq K_{m_j}$, $j\in I_0(\xb)$, are defined as
\begin{equation}\label{def:minface-socp}
	F_j:=\left\{
		\begin{array}{ll}
			K_{m_j} & \textnormal{ if } Dg_j(\xb)d\in \inte(K_{m_j}),\\
			\cone(Dg_j(\xb)d), & \textnormal{ if } Dg_j(\xb)d\in \bd^+(K_{m_j}),\\
			\{0\}, & \textnormal{ if } Dg_j(\xb)d=0.
		\end{array}
	\right.
\end{equation}
and $R\faceq \R^{|I_B(\xb)|}$ is given by
\begin{equation}\label{def:minface-nlp}
	R:=\R^{|I_B(\xb)|}_+ \bigcap_{j\in J} \{c_j\}^\perp, 
\end{equation}
where $c_j$ is the $j$-th vector of the canonical basis of $\R^{|I_B(\xb)|}$, and $J:=\{j\in I_B(\xb)\mid \nabla\phi_j(\xb)^T d=0\}$. Recall that we are assuming for simplicity that $I_B(\xb)=\{1,\ldots,|I_B(\xb)|\}$, and note that $D\G(\xb)d\in F$.

Now, for every $j\in I_0(\xb)$ such that $Dg_j(\xb)d\in \bd^+(K_{m_j})$, let $A_j\in \R^{m_j\times m_j-1}$ be any matrix with full column rank such that $\Im(A_j)=\{Dg_j(\xb)d\}^\perp$, and observe that
\[
	Dg_j(x)^T [F_j^\perp]=\textnormal{span}
	\left(
		\left\{Dg_j(x)^T A_j^i\right\}_{i=1,\ldots,m_j-1}
	\right)
\]
for every such $j$, where $A_j^i$ denotes the $i$-th column of $A_j$. Similarly, for every $j\in I_0(\xb)$ such that $Dg_j(\xb)d=0$, we have 
\[
	Dg_j(x)^T [F_j^\perp]=\textnormal{span}
	(
		\left\{\nabla g_{j,i}(x)\right\}_{i=0,\ldots,m_j-1}
	),
\]
where $\nabla g_{j,i}(x)$ denotes the $i$-th column of $Dg_j(x)\adj$. And for every $j$ such that $Dg_j(\xb)d\in \inte(K_{m_j})$, we have $Dg_j(x)\adj[F_j^\perp]=\{0\}$. Finally, observe that $R^\perp=\spn(\{c_j\}_{j\in J})$ and then
\[
	D\phi(x)^T[R^\perp]=\spn\left(\left\{\nabla \phi_j(x)\right\}_{j\in J} \right).
\]
Therefore, for every $x\in \V$, where $\V$ is the neighborhood of $\xb$ given by Definition~\ref{socp_facial_condition}, the linear space
\begin{equation}\label{socp:eq-facelesspace}
	D\G(x)\adj[F^\perp] =\sum_{j\in I_0(\xb)}Dg_{j}(x)^T [F_j^{\perp}] + D\phi(x)^T[R^\perp]
\end{equation}
is generated by the family of vectors:
\begin{equation}\label{socp:eq-crcqwithoutfaces}
	\bigcup_{\substack{j \in I_0(\xb) \\ Dg_j(\xb)d\in \bdp(K_{m_j}) \\ i=1,\ldots,m_j-1}}\left\{Dg_{j}(x)^{T}A^{i}_{j}\right\}
	\bigcup_{\substack{j \in I_0(\xb) \\ Dg_j(\xb)d=0 \\ i=0,\ldots,m_j-1}} \{\nabla g_{j,i}(x)\} 
	\bigcup_{j \in J} \{\nabla \phi_{j}(x)\},
\end{equation}
which implies that the dimension of~\eqref{socp:eq-facelesspace} equals the rank of~\eqref{socp:eq-crcqwithoutfaces}, for every $x\in \V$. Since this dimension remains constant in $\V$, so does the rank of~\eqref{socp:eq-crcqwithoutfaces}. This means we can apply Proposition~\ref{prop1} to the family of functions
\begin{equation}\label{def:zetasocp}
	\zeta_{i,j}(x):= \left\{ 
		\begin{array}{ll}
		\langle A_{j}^{i}, g_j(x)\rangle, & \textnormal{ if } j \in I_0(\xb), \ Dg_j(\xb)d\in \bdp(K_{m_j}), \ i=1,\ldots,m_j-1,\\
		g_{j,i}(x), & \textnormal{ if } j \in I_0(\xb), \ Dg_j(\xb)d=0, \ i=0,\ldots,m_j-1,\\
		\phi_j(x), & \textnormal{ if } j\in J,
		\end{array}
	\right.
\end{equation}
where $g_{j,i}(x)$ denotes the $i$-th entry of $g_j(x)$ for $j\in J$. Then, consider the following subspace:
\begin{eqnarray*}
\mathcal{S} := \left\{ \begin{tabular}{c|ll}
\multirow{3}{*}{$y \in \mathbb{R}^{n}$} & $A_{j}^{T}Dg_{j}(\xb)y = 0,$    & \textnormal{if} $j \in I_0(\xb), \ Dg_j(\xb)d\in \bdp(K_{m_j})$ \\ 
 										& $Dg_{j}(\xb)y = 0$, & \textnormal{if} $j \in I_0(\xb), \ Dg_j(\xb)d=0$ \\
 										& $\nabla \phi_{j}(\xb)^{T}y = 0 $, & \textnormal{if} $j \in J$ ,
\end{tabular} \right\},
\end{eqnarray*}
and note that $d \in \mathcal{S}$, so it follows that there exists a local diffeomorphism $\psi$ for which items $(i), (ii)$ and $(iii)$ of Proposition~\ref{prop1} are satisfied. Now, define the arc $\xi(t)$ by
\begin{displaymath}
\xi(t) := \psi^{-1}(\xb + td),
\end{displaymath}
for $t\in \R$ small enough so that $\xb+td\in V_2$, where $V_2$ is given by Proposition~\ref{prop1}. Then, we obtain that
\begin{displaymath}
\lim_{t \rightarrow 0^{+}}\xi(t) = \xb, \quad \lim_{t \rightarrow 0^{+}}\dfrac{\xi(t)-\xb}{t} = d.
\end{displaymath}
To complete the proof, it suffices to show that $\xi(t)$ remains feasible for every sufficiently small $t\geq 0$, so this is our goal from this point onwards. Proposition~\ref{prop1} tells us that there exists some $\varepsilon>0$ such that $\zeta_{i,j}(\xi(t))=\zeta_{i,j}(\xb)=0$ for every $t\in (-\varepsilon,\varepsilon)$. In terms of $F$, this means that
\[
	\G(\xi(t))\in \spn(F)
\]
for every such $t$, which follows directly from~\eqref{def:zetasocp}. Now, let us analyse each case separately:
\begin{enumerate}
\item For each index $j\in I_0(\xb)$, consider the Taylor expansion of $g_j(\xi(t))$ around $t=0$, given by
\begin{equation}\label{socp:taylorvertex}
	\begin{aligned}
		g_j(\xi(t)) & =g_j(\xi(0))+tDg_j(\xi(0))\xi'(0)+o(t)\\
		& =g_j(\xb)+tDg_j(\xb)D\psi^{-1}(\xb)d+o(t)\\
		& =tDg_j(\xb)d+o(t)
	\end{aligned}
\end{equation}
Then, we split in three sub-cases:
\begin{itemize}
\item If $Dg_j(\xb)d\in \inte(K_{m_j})$, then it follows from~\eqref{socp:taylorvertex} that $g_j(\xi(t))\in K_{m_j}$ for every $t\in[0,\varepsilon)$, shrinking $\varepsilon$ if necessary;
\item If $Dg_j(\xb)d\in \bdp(K_{m_j})$, then $g_j(\xi(t))\in \spn(Dg_j(\xb)d)$ due to~\eqref{def:minface-socp}, and it follows from~\eqref{socp:taylorvertex} that $g(\xi(t))\in \cone(Dg_j(\xb)d)$ for every $t\in[0,\varepsilon)$, taking a smaller $\varepsilon$ if needed;
\item If $Dg_j(\xb)d=0$, then $g(\xi(t))=0$ for every $t\in [0,\varepsilon)$, due to~\eqref{def:minface-socp}.
\end{itemize}
\item Since $\phi(\xi(t))\in R$ for every $t\in [0,\varepsilon)$, for each index $j\in J$, we have $\phi_j(\xi(t))=0$. On the other hand, for each $j\notin J$, consider the Taylor expansion of $\phi_j(\xi(t))$ around $t=0$:
\begin{equation*}\label{socp:taylorbd}
	\begin{aligned}
		\phi_j(\xi(t)) =\phi_j(\xi(0))+t\nabla \phi_j(\xi(0))^T\xi'(0)+o(t) =t\nabla \phi_j(\xb)^T d+o(t),
	\end{aligned}
\end{equation*}
and since $\nabla \phi_j(\xb)^T d>0$ for every $j\notin J$, it also follows that $\phi_j(\xi(t))>0$ for every $t\in [0,\varepsilon)$, taking a smaller $\varepsilon$ if necessary.
\end{enumerate}
Thus, $\G(\xi(t))\in F$ for every $t\in [0,\varepsilon)$, which also implies that $G(\xi(t))\in \K$ for every such $t$, completing the proof.
\jo{\hfill\Halmos \\}

\si{\end{proof}}

A useful information that can be extracted from the proof above is an equivalent characterization of the facial constant rank property (Definition~\ref{socp_facial_condition}) without faces:

\begin{corollary}\label{socp:crcqwithoutfaces}
Let $\xb\in \Omega$. Then, the facial constant rank property holds at $\xb$ if, and only if, there exists a neighborhood $\V$ of $\xb$ such that: for all subsets $J_1,J_2\subseteq I_0(\xb)$, $J_3\subseteq I_B(\xb)$, such that $m_j>1$ for all $j\in J_1$, and for all $\eta_j\in \bd^+(K_{m_j})$,  $j\in J_1$, the rank of the family 
\begin{equation*}
	\bigcup_{\substack{j \in J_1 \\ i=1,\ldots,m_j}}\left\{Dg_{j}(x)^{T}A^{i}_{j}\right\}
	\bigcup_{\substack{j \in J_2 \\ i=0,\ldots,m_j-1}} \{\nabla g_{j,i}(x)\} 
	\bigcup_{j \in J_3} \{\nabla \phi_{j}(x)\}
\end{equation*}
remains the same for all $x\in \V$, where $A_j\in \R^{m_j\times m_j-1}$ can be any matrix with full column rank such that $\Im(A_j)=\{\eta_j\}^\perp$,  for each $j\in J_1$, and $A^i_j$ denotes the $i$-th column of $A_j$.
\end{corollary}

Before proceeding, we will make a short discussion about Theorem~\ref{socp_facial_theorem} and its implications:

\begin{remark}\label{socp:remlinear}
Note that if all constraints are affine, then every feasible point satisfies the facial constant rank property. Then, it follows from Theorem~\ref{socp_crcq_theorem} that $\tangent = \linear$ in this case, for every $\xb\in \Omega$. We highlight this fact because when it is paired with Example~\ref{socp_example1}, two things can be concluded: first, the facial constant rank property alone is not a CQ for (\ref{socp_multifold}); second, the only reason why constraint linearity is not a CQ for NSOCP is that $H(\xb)$ may not be closed. When $H(\xb)$ is closed, facial constant rank is a CQ, and so is constraint linearity. In other words, the above discussion, in view of the minimality of Guignard's CQ, allows us to conclude that the closedness of $H(\xb)$ is the weakest CQ for linear second-order cone programming problems. 
\end{remark}

The discussion of Remark~\ref{socp:remlinear}, together with Theorem~\ref{socp_crcq_theorem}, motivates our extension of \nlpcr{} (and \nlprcr{}) to NSOCP:

\begin{definition}\label{socp_crcq_definition}
Let $\xb$ be a feasible point of~\eqref{socp_multifold} and let $H(\xb)$ be the set defined in (\ref{h(x)}). We say that the \emph{constant rank constraint qualification for NSOCP} (\socpcr{}) holds at $\xb$, if it satisfies the facial constant rank property and, in addition, the set $H(\xb)$ is closed.
\end{definition}

When $m_{1} = m_{2} = \ldots = m_{q} = 1$, problem~\eqref{socp_multifold} reduces to a NLP problem. Moreover, since the faces of $K_1$ are $\{0\}$ and $\R_+$, the facial constant rank property (Definition~\ref{socp_facial_condition}) reduces to \nlpcr{} in this case, and so does Definition~\ref{socp_crcq_definition}. Moreover, as mentioned before, it follows directly from Theorem~\ref{socp_facial_theorem}, that:

\begin{theorem}\label{socp_crcq_theorem}
The \socpcr{} condition of Definition~\ref{socp_crcq_definition} implies Abadie's CQ.
\end{theorem}

Since the nondegeneracy condition for~\eqref{socp_multifold} holds at a given $\xb\in \Omega$ if, and only if, $D\mathcal{G}(\xb)\adj$ is injective, then by continuity of $D\mathcal{G}$, nondegeneracy implies that $D\mathcal{G}(x)\adj$ remains injective for every $x$ close enough to $\xb$. Therefore, it follows that the nondegeneracy condition implies \socpcr{} as in Definition~\ref{socp_crcq_definition}. However, the converse is not true, as it can be seen in the following example:

\begin{example}\label{socp:ex_crcq_not_rob}
Consider the following constraint
\[ g(x):=(x,x) \in K_{2},\]
at the feasible point $\xb=0$. Since $K_{2}$ is a polyhedral cone and $g$ is linear, then \socpcr{} as in Definition~\ref{socp_crcq_definition} holds at $\xb$. 
However, Robinson's CQ is not satisfied at $\xb = 0$, since 
\[
	Dg(\xb)d=d(1,1)\notin\inte(K_2)
\]
for every $d\in\R$.
Consequently, nondegeneracy is not satisfied, either.
\jo{\hfill\Halmos \\}

\end{example}

Observe that Example~\ref{socp:ex_crcq_not_rob} also shows that \socpcr{} does not imply Robinson's CQ. 
Conversely, Robinson's CQ does not imply \socpcr{} either, meaning they are not related, just as it happens with \nlpcr{} and MFCQ in NLP. Let us show this with an example:

\begin{example}
Consider the constraint:
\[
	g(x):=(x_2,x_1^2) \in K_{2}
\]
at the point $\xb=(0,0)$. Robinson's CQ holds at $\xb$, since $d=(0,1)$ satisfies
\[
	g(\xb) + Dg(\xb)d = (1,0) \in \inte(K_{2}).
\]
On the other hand, take the face $F=\{0\}$ and note that
\[
	Dg(x)\adj[F^\perp] = \spn\left(\left\{
		\begin{bmatrix}
			0 \\ 1
		\end{bmatrix},
		\begin{bmatrix}
			2x_1 \\ 0
		\end{bmatrix}			
	\right\}\right)
\] has dimension $2$ for every $x$ such that $x_1\neq 0$, and dimension $1$ at $\xb$.
\jo{\hfill\Halmos \\}

\end{example}

\begin{remark}\label{socp_remark_eq}
To consider~\eqref{socp_multifold} with an equality constraint in the form $h(x)=0$, where $h\colon \R^n\to \R^p$, one should proceed as in Proposition~\ref{nlp:crcq-conic-improved}. That is, consider 
\[
	g(x):=(g_1(x),\ldots,g_q(x),h(x))
\]
and the cone
\[
	\K:=K_{m_1}\times\ldots\times K_{m_q}\times \{0\}^p.
\]
This will lead to an extension of \nlprcr{}. An extension of the original \nlpcr{} condition can be obtained by writing the equality constraint as a pair of inequality constraints in the form $h(x)\in \R^p_+$ and $-h(x)\in\R^p_+$, just as in Remark~\ref{nlp:remcrcq}, then reducing, and applying Definition~\ref{socp_crcq_definition} to the new reduced cone.
\end{remark}

%
%
\subsection{Strong second-order optimality conditions for NSOCP}

In this subsection we will investigate second-order optimality conditions for (\ref{socp_multifold}) under the facial constant rank property; and, consequently, under \socpcr{} as well. Recall that the second-order condition of Definition~\ref{conic:ssoc} can be further specialized to the context of NSOCP by characterizing the sigma-term explicitly. Following Bonnans and Ram{\'i}rez~\cite{BonRam}, we have for any $\xb\in \Omega$ and any of its associate Lagrange multipliers $\bar\lambda:=(\bar{\lambda}_{1},\ldots, \bar{\lambda}_{q}) \in \Lambda(\xb)$, that
\[
	\displaystyle \sigma(d,\xb,\bar{\lambda})= \sum_{j=1}^{q} d^T\mathcal{H}_{j}(\xb,\bar{\lambda}_{j})d
\] 
for every $d\in C(\xb)$, where 
\begin{equation}\label{socp:sigmah}
\mathcal{H}_{j}(\xb, \bar{\lambda}_{j}) 
:= \left\{
\begin{array}{ll}
	-\dfrac{[\bar{\lambda}_{j}]_{0}}{[g_{j}(\xb)]_{0}}Dg_{j}(\xb)^{T}R_{m_{j}}Dg_{j}(\xb), & \textnormal{if } j\in I_B(\xb), \\ 
	& \\
	0, & \textnormal{otherwise.} 
\end{array}
\right.
\end{equation}

With this in mind, we can prove that SOC holds at $(\xb,\bar\lambda)$ under the facial constant rank property by means of analysing the problem along the curve $\xi(t)$ from the proof of Theorem~\ref{socp_facial_theorem}.

\begin{theorem}\label{socp:thmssoc}
Let $\xb$ be a local minimizer of problem (\ref{socp_multifold}) that satisfies the facial constant rank property. Then, for any given Lagrange multiplier $\bar{\lambda}\in \Lambda(\xb)$, the pair $(\xb,\bar{\lambda})$ satisfies SOC as in Definition~\ref{conic:ssoc}; that is,
\begin{equation}\label{socp:eqsoc}
	d^T \nabla^2 f(\xb)d+\sum_{j=1}^q \left\langle D^{2}g_{j}(\xb) [d,d], \bar{\lambda}_{j}\right\rangle - \sigma(d,\xb,\bar\lambda)\geq 0,
\end{equation}
for every $d\in C(\xb)=\linear\cap \{\nabla f(\xb)\}^\perp$.
\end{theorem}

\jo{\proof{Proof.}}\si{\begin{proof}}
If $\Lambda(\xb)=\emptyset$, the result holds trivially. Otherwise, let $\bar{\lambda}:=(\bar{\lambda}_{1}, \ldots, \bar{\lambda}_{q})\in \Lambda(\xb)$ be arbitrary and fixed. Our aim is to prove that inequality~\eqref{conic:eqsoc} holds for the pair $(\xb,\bar{\lambda})$, for every $d\in C(\xb)$. So let $d \in C(\xb)$ be also arbitrary, and let $F$ be as in~\eqref{def:minface-full}. Recall that, for the sake of simplicity and without loss of generality, we are assuming $I_B(\xb)=\{1,\ldots,|I_B(\xb)|\}$.

Proceeding in the same way as in the proof of Theorem~\ref{socp_facial_theorem}, since the facial constant rank property holds at $\xb$ and $d\in \linear$, we can construct a curve $\xi:(-\varepsilon,\varepsilon) \rightarrow \R^{n}$, for some $\varepsilon>0$, such that: $\xi(0) = \xb$, $\xi^{\prime}(0) = d$, and $\G(\xi(t))\in \spn(F)$ for every $t\in (-\varepsilon,\varepsilon)$. In addition, $\G(\xi(t))\in F$ for every $t\in[0,\varepsilon)$, meaning $\xi(t)$ is feasible for all such $t$. Since $\xb$ is a local minimizer of~\eqref{socp_multifold}, then $t=0$ is a local minimizer of the function $\varphi(t) := f(\xi(t))$ subject to the constraint $t\geq 0$. Then, it is easy to see that
\begin{equation}\label{2a_arc}
\varphi^{\prime \prime}(0) = d^{T} \nabla^{2}f(\xb)d + \nabla f(\xb)^{T} \xi^{\prime \prime}(0) \geq 0.
\end{equation}

The rest of this proof consists of computing $\nabla f(\xb)^{T} \xi^{\prime \prime}(0)$. To do this, we will use an auxiliary complementarity function defined as
\begin{displaymath}
R(t) := \sum_{j \in I_{0}(\xb)}\langle g_{j}(\xi(t)), \bar{\lambda}_{j} \rangle + \sum_{j \in I_{B}(\xb)} [\bar{\lambda}_{j}]_{0}\phi_j(\xi(t)).
\end{displaymath}
First, we claim that $R(t)=0$ for every $t\in (-\varepsilon,\varepsilon)$. Let us prove this, analysing each term separately:
\begin{enumerate}
\item For each $j\in I_0(\xb)$, it follows from the KKT conditions that
\begin{equation}\label{eq:compcrit-socp}
	\langle Dg_j(\xb)d, \bar{\lambda}_j \rangle =  \langle d, Dg_j(\xb)^T\bar{\lambda}_j \rangle = \langle d, -\nabla f(\xb)\rangle = 0,
\end{equation}
which implies the following:
\begin{itemize}
\item If $Dg_j(\xb)d\in \inte(K_{m_j})$, then $\bar{\lambda}_j=0$, since $\bar{\lambda}_j\in K_{m_j}\pol$;
\item If $Dg_j(\xb)d\in \bdp(K_{m_j})$ we have $g_j(\xi(t))\in \spn(Dg_j(\xb)d)$, and consequently, $\langle g_j(\xi(t)),\bar{\lambda}_j \rangle=0$ for every $t\in(-\varepsilon,\varepsilon)$ due to~\eqref{eq:compcrit-socp};
\item If $Dg_j(\xb)d=0$, then $g(\xi(t))=0$ also for every $t\in (-\varepsilon,\varepsilon)$, due to~\eqref{def:minface-socp}.
\end{itemize}
The above reasoning implies that $\langle g_{j}(\xi(t)), \bar{\lambda}_{j} \rangle=0$ for every $t\in [0,\varepsilon)$ and every $j\in I_0(\xb)$.  
\item For each $j\in I_B(\xb)$, consider $J$ as in~\eqref{def:minface-nlp} and it follows that if $\nabla \phi_j(\xb)^Td=0$, then $\phi_j(\xi(t))=0$ for every $t\in (-\varepsilon,\varepsilon)$. On the other hand, recall that $\bar{\lambda}_j=\frac{[\bar{\lambda}_{j}]_{0}}{[g_{j}(\xb)]_{0}}R_{m_j} g_j(\xb)$ for every $j\in I_B(\xb)$ due to complementarity, but using~\eqref{socp:multformat}, we see that
\[
	[\bar{\lambda}_{j}]_{0}\langle \nabla \phi_j(\xb), d \rangle=\frac{[\bar\lambda_j]_0}{[g_j(\xb)]_0}\langle Dg_j(\xb)^T R_{m_j}g_j(\xb), d \rangle = \langle \bar{\lambda}_j, Dg_j(\xb)d \rangle=0.
\]
Therefore, if $\langle \nabla \phi_j(\xb), d \rangle>0$, then $[\bar{\lambda}_j]_0=0$.
\end{enumerate}

Knowing that $R(t)=0$ for every $t\in (-\varepsilon,\varepsilon)$, we obtain that the derivatives of $R(t)$ are also zero for all such $t$. Let us compute them: the first derivative of $R(t)$ is given by

\begin{eqnarray*}
R^{\prime}(t) &=& \sum_{j \in I_{0}(\xb)} \left\langle Dg_{j}(\xi(t)) \xi^{\prime}(t), \bar{\lambda}_{j}\right\rangle + \sum_{j \in I_{B}(\xb)} [\bar{\lambda}_{j}]_{0} \left\langle \nabla \phi_j(\xi(t)),\xi^{\prime}(t)\right\rangle,
\end{eqnarray*}
and the derivative of $R'(t)$ is
\begin{eqnarray*}
R^{\prime \prime}(t) &=& \sum_{j \in I_{0}(\xb)} \left\langle D^{2}g_{j}(\xi(t)) [\xi^{\prime}(t), \xi^{\prime}(t)], \bar{\lambda}_{j}\right\rangle 
+
\sum_{j \in I_{0}(\xb)} \left\langle Dg_{j}(\xi(t))\adj\bar{\lambda}_{j}, \xi^{\prime \prime}(t)\right\rangle\\
& & +
\sum_{j \in I_{B}(\xb)} [\bar{\lambda}_{j}]_{0} \left(\left\langle D^2 \phi_j(\xi(t))\xi'(t),\xi^{\prime}(t)\right\rangle
+\left\langle \nabla \phi_j(\xi(t)),\xi^{\prime\prime}(t)\right\rangle\right).
\end{eqnarray*}
Since $R^{\prime\prime}(t)$ is continuous, taking the limit $t\to 0$, we obtain
\begin{eqnarray*}
R^{\prime \prime}(0)=\lim_{t\to 0^+} R^{\prime\prime}(t) 
&=& 
\sum_{j \in I_{0}(\xb)} \left\langle D^{2}g_{j}(\xb) [d,d], \bar{\lambda}_{j}\right\rangle  
+
\sum_{j \in I_{0}(\xb)} \left\langle Dg_{j}(\xb)\adj\bar{\lambda}_{j}, \xi^{\prime \prime}(0)\right\rangle \\
& &+
\sum_{j \in I_{B}(\xb)} [\bar{\lambda}_{j}]_{0} \left(\left\langle D^2 \phi_j(\xb)d,d\right\rangle
+\frac{1}{[g_j(\xb)]_0}\left\langle Dg_j(\xb)^T R_{m_j}g_j(\xb),\xi^{\prime\prime}(0)\right\rangle\right).
\end{eqnarray*}
The above expression can be simplified using the relation
\[
	\begin{aligned}
		\left\langle D^2 \phi_j(\xb)d,d\right\rangle 
		&= \frac{\langle \widehat{Dg_j(\xb)d},\widehat{g_j(\xb)} \rangle^2}{\|\widehat{g_j(\xb)}\|^3}-\frac{\|\widehat{Dg_j(\xb)d}\|^2}{\|\widehat{g_j(\xb)}\|} +\left\langle D^2g_j(\xb)[d,d], \nabla \phi_j(\xb) \right\rangle\\		
		&=\frac{1}{[g_j(\xb)]_0}\left(\left\langle R_{m_j}Dg_j(\xb)d, \ Dg_j(\xb)d\right\rangle +\left\langle D^2g_j(\xb)[d,d], R_{m_j}g_j(\xb) \right\rangle\right),
	\end{aligned}
\]
which can be directly computed from the definition of $\phi_j$, $j\in I_B(\xb)$, since in this case $[g_j(\xb)]_0=\|\widehat{g_j(\xb)}\|$ and $\langle Dg_j(\xb)d, R_{m_j}g_j(\xb)\rangle=0$. Then, we get
\begin{eqnarray}\label{new-Rd2}
R^{\prime \prime}(0) &=& \sum_{j \in I_{0}(\xb)\cup I_B(\xb)} \left\langle D^{2}g_{j}(\xb) [d,d], \bar{\lambda}_{j}\right\rangle + \sum_{j \in I_{0}(\xb)\cup I_B(\xb)} \left\langle Dg_{j}(\xb)^{T}\bar{\lambda}_{j},\xi^{\prime \prime}(0)\right\rangle \nonumber \\
 & & + \sum_{j \in I_B(\xb)}\dfrac{[\bar{\lambda}_{j}]_{0}}{[g_{j}(\xb)]_{0}} \left\langle R_{m_j}Dg_j(\xb)d, \ Dg_j(\xb)d\right\rangle = 0.
\end{eqnarray}
Moreover, by the KKT conditions, we have
\[
	\nabla f(\xb)^{T} \xi^{\prime \prime}(0) = - \sum_{j \in I_{0}(\xb)\cup I_B(\xb)} \left\langle Dg_{j}(\xb)^{T}\bar{\lambda}_{j},\xi^{\prime \prime}(0)\right\rangle ,
\]
which yields together with equation~\eqref{new-Rd2}, the following:
\begin{equation}\label{socp:eqterm}
	\nabla f(\xb)^{T} \xi^{\prime \prime}(0) = \sum_{j \in I_{0}(\xb)\cup I_B(\xb)} \left\langle D^{2}g_{j}(\xb) [d,d], \bar{\lambda}_{j}\right\rangle + \sum_{j \in I_{B}(\xb)}\dfrac{[\bar{\lambda}_{j}]_{0}}{[g_{j}(\xb)]_{0}} d^{T}Dg_{j}(\xb)^{T}R_{m_{j}}Dg_{j}(\xb)d.
\end{equation}
Therefore, since $\bar\lambda_j=0$ for every $j\in I_{\inte}(\xb),$
\begin{displaymath}
d^{T} \nabla^{2}f(\xb)d +\sum_{j=1}^q \left\langle D^{2}g_{j}(\xb) [d,d], \bar{\lambda}_{j}\right\rangle -\sigma(d,\xb,\bar{\lambda}) \geq 0. 
\end{displaymath}
Since $d \in C(\xb)$ is arbitrary, we conclude that $\xb$ satisfies SOC with respect to $\bar\lambda$, which was also chosen arbitrarily and remained fixed from the very beginning. Thus, the proof is complete.
\jo{\hfill\Halmos \\}

\si{\end{proof}}

Observe that Theorem~\ref{socp:thmssoc} implies that the facial constant rank property ensures the fulfilment of the strong second-order necessary condition at a given point $\xb$, in the sense that for every $\bar{\lambda}\in \Lambda(\xb)$, and every $d\in C(\xb)$, inequality~\eqref{conic:eqsoc} holds true. If, in addition, $H(\xb)$ is closed (CRCQ), then $\Lambda(\xb)\neq \emptyset$, and as consequence, we obtain that the strong second-order condition is satisfied in the presence of \socpcr{}. It is also worth mentioning that since the strong necessary condition of Theorem~\ref{socp:thmssoc} implies the classical condition of~Theorem~\ref{conic:bsoc}, then it also induces a sufficient second-order optimality condition after replacing $\geq$ by $>$ in inequality~\eqref{socp:eqsoc}.

\begin{remark}\label{socp:remarkzz}
In contrast with the facial constant rank property, the condition presented in~\cite[Definition 2.1]{ZZ} fails to be a CQ even when $H(\xb)$ is closed. In fact, let us recall the counterexample presented in~\cite{andreani2020erratum}:
\begin{equation}
\begin{array}{ll}
\nonumber\mbox{\textnormal{Minimize }} 	& f(x):=-x,\\ 
\mbox{s.t. }	& g(x):=(x,x+x^2) \in K_{2},
\end{array}
\end{equation}
The unique solution of this problem is $\xb = 0$. For this particular example,~\cite[Definition 2.1]{ZZ} holds if, and only if, $\{1, 1+2x\}$ remain with constant rank in some neighborhood of $\xb$ (one may consider also all of its subfamilies, see \cite{andreani2020erratum}). Of course, this is verified, and since $K_2$ is polyhedral, the set $H(\xb)$ is closed. However, $\xb$ does not satisfy the KKT conditions. 

On the other hand, to see that  \socpcr{} as in Definition~\ref{socp_crcq_definition} does not hold at $\xb$, take $F:=\cone((1,1))\faceq K_2$ and note that
\begin{displaymath}
Dg(x)^{T}[F^\perp] =\spn(-2x)
\end{displaymath}
has dimension $1$ for every $x\neq 0$, but has dimension zero at $\xb$. In particular, this example shows that~\socpcr{} as in Definition~\ref{socp_crcq_definition} is not a mere correction of the condition presented in~\cite{ZZ}, and that the condition of~\cite{ZZ} cannot be corrected by simply adding the closedness of $H(\xb)$ to its definition.
\end{remark}

\subsection{About the sequential constant rank CQ}

In~\cite{seqcrcq-socp}, we introduced an alternative extension of \nlpcr{} for (\ref{socp_multifold}) that was based on a special re-characterization of the nondegeneracy condition \cite{weak-sparse-cq} in terms of the eigenvectors of some perturbations of $g(\xb)$. Let us recall an equivalent characterization of it, which will be used here as a definition for simplicity.

\begin{definition}[Seq-CRCQ for NSOCP]\label{socp:seqcrcq}
Let $\xb\in \Omega$. We say that the \emph{Sequential-CRCQ} (Seq-CRCQ) condition holds at $\xb$ if for every vector $\bar{w}_j\in \R^{m_j-1}$ with $\|\bar{w}_j\|=1$, $j\in I_0(\xb)$, there is a neighborhood $\V$ of $(\xb,\bar{w})$, $\bar{w}:=(\bar{w}_j)_{j\in I_0(\xb)}$, such that: for all subsets $J_1,J_2\subseteq I_0(\xb)$ and $J_3\subseteq I_B(\xb)$, if the family
\[
	\mathcal{D}(x,w):=
	\left\{ Dg_j(x)^T \left(1, -{w}_j\right) \right\}_{j\in J_1}
	\bigcup
	\left\{ Dg_j(x)^T \left(1, {w}_j\right) \right\}_{j\in J_2}
	\bigcup
	\left\{ Dg_j(x)^T \left(1, -\frac{\widehat{g_j(x)}}{\|\widehat{g_j(x)}\|} \right) \right\}_{j\in J_3}
\] is linearly dependent at $(x,w):=(\xb,\bar{w})$, then $\mathcal{D}(x,w)$ remains linearly dependent for all $(x,w)\in \V$ such that $\|w_j\|=1$, $j\in J_1\cup J_2$, where $w:=(w_j)_{j\in I_0(\xb)}$.
\end{definition}

This constraint qualification was used in~\cite{seqcrcq-socp} to achieve global convergence of a class of algorithms to KKT points, and some interesting properties were shown together with a weaker variant of Seq-CRCQ. Namely, it is also independent of Robinson's CQ, strictly weaker than nondegeneracy, and it implies the \emph{metric subregularity CQ} (also known as \emph{error bound CQ}). Moreover, note that if $I_0(\xb)=\emptyset$, then Seq-CRCQ coincides with the facial constant rank property, which in turn coincides with \socpcr{}. However, this is not necessarily true otherwise. In the following example, we show that \socpcr{} according to Definition~\ref{socp_crcq_definition} does not imply Seq-CRCQ.

\begin{example}\label{socp:exseq}
Consider the constraint:
\begin{equation*}\nonumber
g(x) = (x, -x, 0)\in K_{3},
\end{equation*}
and let $\xb = 0$, a feasible point. Since $g$ is affine, then the facial constant rank property holds at $\xb$ (see Remark~\ref{socp:remlinear}). Now, let us show that $H(\xb)$ is closed: since $g(\xb) = 0$, it holds that
\[
	H(\xb) = Dg(\xb)^{T} K_{3} = \{v_{1}-v_{2} \mid (v_{1}, v_{2}, v_{3}) \in K_{3}\}=\R_+.
\]
Therefore, $H(\xb)$ is a closed set, and \socpcr{} according to Definition~\ref{socp_crcq_definition} holds at $\xb$.

On the other hand, Seq-CRCQ does not hold at $\xb$, since for any $w=(w_1,w_2)\in \R^2$,
\[
	Dg(\bar{x})^T (1,w)=1-w_1 \quad \textnormal{and} \quad Dg(\bar{x})^T (1,-w)=1+w_1;
\]
then, take $\bar{w}=(1,0)$ and any sequence $\seq{w^k}\to \bar{w}$ such that $w_1^k\neq 1$ for all $k\in \N$, to see that $Dg(\bar{x})^T (1,w^k_1)\neq 0$ for every $k\in \N$, but $Dg(\bar{x})^T (1,\bar{w})=0$.
\jo{\hfill\Halmos \\}

\end{example}

We were not able to prove nor find a counterexample for the converse statement. However, with only Example~\ref{socp:exseq} at hand, we already know that~\socpcr{} is in the worst case independent of Seq-CRCQ, and in the best case, strictly weaker than it, meaning the results of this paper either improve or are parallel to the results of~\cite{seqcrcq-socp}.

%
%

\section{Nonlinear semidefinite programming}\label{sec:sdp}

In this section, we will study constant rank conditions for nonlinear semidefinite programming problems, which can be stated in standard form as follows:
\begin{equation}\label{nsdp}
\begin{array}{ll}
\nonumber\mbox{\textnormal{Minimize }} & f(x),\\ 
\mbox{s.t. }	& G(x)\succeq 0.
\end{array}
\tag{NSDP}
\end{equation}
This problem can be seen as a particular case of~\eqref{conic_problem}, letting $\E=\S^m$ be the space of $m\times m$ symmetric matrices with real entries, and 
\[
	\K=\mathbb{S}^m_+:=\{A\in \S^m\mid z^T Az\geq 0, \ \forall z\in \R^m\}
\]
be the cone of all $m\times m$ symmetric positive semidefinite matrices, with $G:\R^n\to\mathbb{E}$ being twice continuously differentiable. The symbol $\succeq$ denotes the partial order induced by $\S^m_+,$ meaning that $A\succeq B$ if, and only if, $A-B\in \mathbb{S}^m_+$. In this section, for any given $A\in \S^m$ we will denote by $\mu_i(A)$ the $i$-th eigenvalue of $A$ arranged in non-increasing order, and $u_i(A)$ will denote an associated unitary eigenvector.

Recall from Section~\ref{sec:nlp} that the constant rank constraint qualification can be obtained in two steps: first, reduce the problem to consider only the locally relevant part of the constraint; then, analyse the image of the faces of the reduced cone by the derivative of the reduced constraint function. Therefore, we begin by recalling a matrix analysis lemma that will be useful for the reduction part.

\begin{lemma}
\label{lemma:convergence}
Let $\bar{A}\succeq 0$ and denote by $r$ the rank of $\bar{A}$. Also, let $\bar{E} \in \R^{m\times m-r}$ be a matrix whose columns form an orthonormal basis of $\Ker(\bar{A})$. Then, there exists an analytic matrix function $\mathcal{E}\colon  \mathbb{S}^m\to \R^{m\times m-r}$ such that $\mathcal{E}(\bar{A})=\bar{E}$ and, for all $A$ close enough to $\bar{A}$, the columns of $\mathcal{E}(A)$ form an orthonormal basis of the space spanned by the eigenvectors associated with the $m-r$ smallest eigenvalues of $A$.
\end{lemma}

\jo{\proof{Proof.}}\si{\begin{proof}}
Although this proof can be found in~\cite[Example 3.140]{bonnans-shapiro} and~\cite[Section 2.3]{auslender}, we shall include it here for completeness purposes. For any given $A\in \S^m$, let $S(A)$ be the space spanned by the eigenvectors $u_{r+1}(A),...,u_{m}(A)$, and let $\Pi(A)$ denote the orthogonal projection matrix onto $S(A)$. It is known that $\Pi(A)$ is an analytic function of $A$ in a sufficiently small neighborhood of $\bar{A}$ (see, for example, \cite[Theorem 1.8]{kato}). Consequently, the function $V(A):=\Pi(A)\bar{E}$ is also an analytic function of $A$ in a neighborhood of $\bar{A}$, and moreover $V(\bar{A})=\bar{E}$.
  It follows that for all $A$ sufficiently close to $\bar{A}$, the rank of
  $V(A)$ is equal to the rank of $V(\bar{A})=\bar{E}$, so the $m-r$ columns of $V(A)$ are linearly independent, and $\Im(V(A))=S(A)$, when $A$ is sufficiently close to $\bar{A}$.
  Now, let $\mathcal{E}(A)$ be a matrix whose columns are obtained by applying
  the Gram-Schmidt orthonormalization process to the columns of
  $V(A)$. The matrix $\mathcal{E}(A)$ is well-defined and also analytic
   in a neighborhood of $\bar{A}$. Moreover, it satisfies $\mathcal{E}(A)^T \mathcal{E}(A)= \I_{m-r}$ and $\Im(\mathcal{E}(A)) = S(A)$, for all $A$ sufficiently close to $\bar{A}$, and $\mathcal{E}(\bar{A})=V(\bar{A})=\bar{E}$, which concludes the proof. Notice however that the columns of $\mathcal{E}(A)$ are not necessarily eigenvectors of $A$.
\jo{\hfill\Halmos \\}

\si{\end{proof}}
 
Now, let $\xb\in \Omega$, denote by $r$ the rank of $G(\xb)$, and let $\bar{E}\in \R^{m\times m-r}$ be an arbitrary matrix with orthonormal columns that span $\Ker (G(\xb))$. Moreover, let $\mathcal{E}$ be the analytic matrix function with the properties described in Lemma~\ref{lemma:convergence}, such that $\mathcal{E}(G(\xb))=\bar{E}$. Observe that $\S^m_+$ is reducible to 
\[
	\mathcal{C}:=\S^{m-r}_+
\] in a neighborhood $\viz$ of $G(\xb)$ by the mapping $\Xi\colon \viz\to \S^{m-r}$ given by
\[
	\Xi(Y):=\mathcal{E}(Y)^T Y \mathcal{E}(Y),
\]
for every $Y\in\viz$ close enough to $G(\xb)$ so that $\mu_i(Y)>0$ for every $i=1,\ldots,r$. Then, define the function $E:=\mathcal{E}\circ G$, consider the reduced constraint function
\[
\mathcal{G}(x):= E(x)^T G(x)E(x),
\]
and for every $x$ sufficiently close to $\xb$, we have that $G(x)\in \S^m_+$ if, and only if, $\mathcal{G}(x)\in \S^{m-r}_+$. Moreover, it is worth recalling that, since the function $\mathcal{E}$ is analytic, the degree of differentiability of $\mathcal{G}$ is the same as of $G$.

Following Bonnans and Shapiro~\cite[Equation 5.161]{bonnans-shapiro}, we see that the linearized cone of the original constraints of~\eqref{nsdp} at $\xb\in \Omega$ can be written as
\[
	\linear = \left\{ d\in \mathbb{R}^n \mid \bar{E}^T DG(\xb)d\bar{E}\succeq 0 \right\},
\]
which also coincides with the linearized cone of the reduced constraint at $\xb$, because $E(\xb)=\bar{E}$ and for each $x$ close enough to $\xb$, we have
\begin{align*}
	D\mathcal{G}(x)[ \ \cdot \ ] & = DE(x)[ \ \cdot \ ]^T G(x)E(x)+E(x)^T  DG(x)[ \ \cdot \ ]E(x) +  E(x)^T G(x)DE(x)[ \ \cdot \ ],
\end{align*}
so $D\mathcal{G}(\xb)[ \ \cdot \ ]=\bar{E}^T DG(\xb)[ \ \cdot \ ]\bar{E}$. For more details on this reduction approach, see~\cite{bonn-comi-shap,myreportwithbonnans}.

In the next section, we will introduce a constant rank-type condition for NSDP by means of the faces of the reduced cone. 

%
%
\subsection{A facial constant rank constraint qualification for NSDP}\label{sec:crcqsdp}

Following the exposition of Pataki~\cite{pataki-geometry}, the faces of $\mathcal{C}=\S^{m-r}_+$ can be represented in a very simple way: $F$ is a face of $\S^{m-r}_+$ if, and only if, there exists an orthogonal matrix $U\in \R^{{m-r}\times {m-r}}$ and some $s\in \{1,\ldots,{m-r}\}$ such that
\[
	F= 
	\left\{
		U
		\begin{bmatrix}
			A_{11} & 0 \\ 0 & 0
		\end{bmatrix}
		U^T
		\ \middle| \ A_{11}\in \S^s_+	
	\right\}.
\]

With this in mind, let us define the analogue of~Definition~\ref{socp_crcq_definition} for NSDP:

\begin{definition}\label{sdp_crcqdef}
Let $\xb\in\Omega$ and let $r$ be the rank of $G(\xb)$. We say that the \emph{facial constant rank} property holds at $\xb$ if there exists a neighborhood $\V$ of $\xb$ such that: for each $F\faceq\S^{m-r}_+$, the dimension of $D\mathcal{G}(x)\adj [F^\perp]$ remains constant for every $x\in \V$.
\end{definition}

Following the discussion after Proposition~\ref{nlp:crcq-conic-improved} and also after Definition~\ref{socp_facial_condition}, to better visualize the meaning of Definition~\ref{sdp_crcqdef}, recall that $\Omega$ is locally equivalent to $\G^{-1}(\S^{m-r}_+)$ and that the faces of $\S^{m-r}_+$ can be regarded as linear approximations of $\S^{m-r}_+$, in some sense. Then, for every $F\faceq \S^{m-r}$, the set $D\G(x)^{-1}(\spn(F))$ defines a possible linear approximation of $\Omega$ around $\xb$. The reasoning after Proposition~\ref{nlp:crcq-conic-improved} still holds in the context of NSDP and it follows that the facial constant rank property holds at $\xb\in \Omega$ if, and only if, the dimension of $D\G(x)^{-1}(\spn(F))$ remains constant for all $x$ in a neighborhood of $\xb$, at every $F\faceq\S^{m-r}_+$. From this point of view, the facial constant rank property demands all linear approximations of the feasible set to remain with constant dimension in the vicinity of $\xb$. 

Now, we proceed to the main result of this section.

\begin{theorem}\label{sdp:fcrtl}
Let $\xb\in \Omega$. If $\xb$ satisfies the facial constant rank property, then $\tangent=\linear$.
\end{theorem}

\jo{\proof{Proof.}}\si{\begin{proof}}
Let $d\in \linear$, denote by $s$ the rank of $\bar{E}^T DG(\xb)d\bar{E}$, and let $\bar{Q}\in \R^{m-r\times m-r}$ be an orthogonal matrix such that 
\[
	\bar{Q}^T \bar{E}^T DG(\xb)d \bar{E} \bar{Q} = \begin{bmatrix}
		R & 0\\ 0 & 0
	\end{bmatrix},
\] where $R\succ 0$ is an $s\times s$ diagonal matrix. Let $\bar{W}$ be the matrix formed by the columns of $\bar{Q}$ corresponding to the positive eigenvalues of $\bar E^T DG(\xb)d\bar E$; that is, $\bar{W}^T \bar{E}^T DG(\xb)d\bar{E} \bar{W} = R$. Then, consider the face of $\S^{m-r}_+$ given by
\begin{equation}\label{sdp:faceproof}
	F:=
	\left\{
		\bar{Q}
		\begin{bmatrix}
			A_{11} & 0 \\ 0 & 0
		\end{bmatrix}
		\bar{Q}^T
		\ \middle| \ A_{11}\in \S^s_+	
	\right\}
\end{equation}
and note that $\bar{E}^T DG(\xb)d\bar{E}\in F$.
Let $\eta_1,\ldots,\eta_N$ be a basis of $F^\perp$, where $N:=\dim(F^\perp)$, and note that
\begin{equation}\label{eq:rank}
	D\mathcal{G}(x)\adj [F^\perp]=
	\spn\left(\left\{
		D\mathcal{G}(x)^T[\eta_i]\right\}_{i\in \{1,\ldots,N\}}\right).
\end{equation}
Therefore, the facial constant rank property can be equivalently stated as the constant rank of the family 
\[\left\{
		D\mathcal{G}(x)^T[\eta_i]\right\}_{i\in \{1,\ldots,N\}}\]
in a neighborhood of $\xb$. Furthermore, let $\zeta_i(x):=\langle \mathcal{G}(x),\eta_i\rangle$ and note that
\[
	\nabla \zeta_i(x)=D\mathcal{G}(x)\adj[\eta_i]
\] for every $i\in \{1,\ldots,N\}$.

 Then, by Proposition~\ref{prop1}, there exist two neighborhoods $V_1$ and $V_2$ of $\xb$, and a curve $\psi\colon V_1\to V_2$ such that $\psi(\xb)=\xb$, $D\psi(\xb)=\I_n$, and $\zeta_i(\psi^{-1}(\xb+y))=\zeta_i(\xb)$ for every $i\in \{1,\ldots,N\}$ and every $y$ in the subspace
\[
	\mathcal{S}:=\left\{
		y\in \R^n\mid \langle \nabla \zeta_i(\xb), y\rangle=0, \ \forall i\in \{1,\ldots,N\}
	\right\}.
\]
Since $D\mathcal{G}(\xb)d\in F$, we see that $\langle d, D\mathcal{G}(\xb)\adj[\eta_i]\rangle=\langle D\mathcal{G}(\xb)d,\eta_i\rangle=0$ for every $i\in\{1,\ldots,N\}$, so $d\in \mathcal{S}.$ Then, let $\varepsilon>0$ be such that $\xb+td\in V_2$ for every $t\in (-\varepsilon,\varepsilon)$, and define $\xi(t):=\psi^{-1}(\xb+td)$ for every such $t$. Moreover, note that $\xi'(t)=d$ and $\xi(0)=\xb$. 

Now, for every $t\in (-\varepsilon,\varepsilon)$, we have
$\langle\mathcal{G}(\xi(t)),\eta_i\rangle=0$ for every $i\in \{1,\ldots,N\}$, whence follows that 
\[
	\mathcal{G}(\xi(t))\in \spn(F)
\] 
for every such $t$, meaning also
\[
	\bar{Q}^T \mathcal{G}(\xi(t))\bar{Q} = \begin{bmatrix}
		\bar{W}^T \mathcal{G}(\xi(t))\bar{W} & 0\\
		0 & 0
	\end{bmatrix}.
\]
On the other hand, considering the Taylor expansion of $\mathcal{G}(\xi(t))$ around $t=0$,
\[
	\mathcal{G}(\xi(t))=\bar{E}^T G(\xb)\bar{E} + t \bar{E}^T DG(\xb)d \bar{E}+ o(t)= t \bar{E}^T DG(\xb)d \bar{E}+ o(t),
\]
we observe that
\[
	\bar{W}^T \mathcal{G}(\xi(t))\bar{W}= t \bar{W}^T\bar{E}^T DG(\xb)d \bar{E}\bar{W} + o(t)=tR+o(t)\succ 0,
\]
for $t\in[0,\varepsilon)$, shrinking $\varepsilon$ if necessary. Thus, 
\[
	\mathcal{G}(\xi(t))\in F\subseteq \S^{m-r}_+
\]
for every $t\in[0,\varepsilon)$, and then $G(\xi(t))\succeq 0$ for all such $t$. Therefore, it follows that $d\in \tangent$.
\jo{\hfill\Halmos \\}

\si{\end{proof}}

\begin{remark}\label{sdp:remlinear}
Similarly to Remark~\ref{socp:remlinear}, we observe that if $G$ is affine, then every $\xb\in \Omega$ satisfies the facial constant rank property, which implies it is not a CQ on its own, unless $H(\xb)$ as defined in \eqref{h(x)} is closed. We remark this fact because it implies that the weakest CQ that guarantees zero duality gap in linear SDP problems is the closedness of $H(\xb)$. 
\end{remark}

With this in mind, we present our extension of \nlpcr{} (and \nlprcr{}) for NSDP:

\begin{definition}[\sdpcr{}]\label{sdp:crcq_def}
Let $\xb\in \Omega$. We say that $\xb$ satisfies the \emph{constant rank constraint qualification condition for NSDP} (\sdpcr{}) if it satisfies the facial constant rank property and, in addition, $H(\xb)$ is closed.
\end{definition}

And, as an immediate consequence of Theorem~\ref{sdp:fcrtl}, we obtain the following:

\begin{theorem}
Let $\xb\in \Omega.$ If $\xb$ satisfies \sdpcr{}, then it also satisfies Abadie's CQ.
\end{theorem}

Now, we are led to compare our \sdpcr{} condition with other CQs from the literature. First, let us show that it is, in general, independent of Robinson's CQ.

\begin{example}\label{sdp:exseq}
Consider the following constraint:
\[
	G(x) := \left[ \begin{array}{cc} -x & 0 \\ 0 & x \end{array} \right].
\]
The only feasible point is $\xb = 0$, for which one has $G(\xb) = 0 \in \mathbb{S}^{2}_{+}$. Since $G$ is linear it is enough to show that $H(\xb)$ is closed (see Remark~\ref{sdp:remlinear}).  In this case, 
\[
	H(\xb) = \left\{DG(\xb)\adj A \ \middle| \ A=\begin{bmatrix}
	a_{11} & a_{12} \\ a_{12} & a_{22}
\end{bmatrix} \in \mathbb{S}^{2}_{+}\right\}.
\] Since $DG(\xb)\adj A = \langle DG(\xb), A \rangle = a_{22}-a_{11}$, but $a_{11}$ and $a_{22}$ are both nonnegative, it follows that $H(\xb)=\R$, which is closed. On the other hand, Robinson's CQ does not hold at $\xb$. In fact, given a real number $d \in \R$, we have that 
\[
	G(\xb) + DG(\xb)d = \left[ \begin{array}{cc} -d & 0 \\ 0 & d \end{array} \right],
\]
which is not in $\textnormal{\inte}(\mathbb{S}^{2}_{+})$ regardless of $d \in \R$.
\jo{\hfill\Halmos \\}

\end{example}

The above example shows that \sdpcr{} does not imply Robinson's CQ. Conversely, we will show in the next example, that Robinson's CQ does not imply \sdpcr{} either.

\begin{example}\label{sdp:exreduction}
Consider the following constraint given by
\[
G(x) := \left[ \begin{array}{cc} x_2 & x_1^{2} \\ x_1^{2} & x_2 \end{array} \right]
\]
at the point $\xb = (0,0)$. Then, for any direction $d = (d_{1},d_{2}) \in \R^{2}$, it follows that 
\[
	DG(x)d = \left[ \begin{array}{cc} 0 & 2x_1 \\ 2x_1 & 0 \end{array} \right] d_{1} + \left[ \begin{array}{cc} 1 & 0 \\ 0 & 1 \end{array} \right] d_{2}.
\]

Now, consider $\overline{d}=(0,1)$. Then
\[
	G(\xb) + DG(\xb)\overline{d} = \left[ \begin{array}{cc} 1 & 0 \\ 0 & 1 \end{array} \right] \in \textnormal{\inte}(\mathbb{S}^{2}_{+})
\]
and Robinson's CQ holds at $\xb$. In order to analyze whether \sdpcr{} holds at $\xb$ or not, let $\bar{E}\in\R^{2\times 2}$ be a given orthogonal matrix and let us define $\mathcal{G}(x):=
\bar{E}^T G(x)\bar{E}$. Then, take
\[
	F:=\left\{
		\bar{E}^T
		\begin{bmatrix}
		a & 0 \\ 0 & 0
		\end{bmatrix}
		\bar{E}
		\ \middle| \
		a\geq 0
	\right\}\faceq \S^2_+
\]
and note that
\[
	D\G(x)\adj[F^\perp]=DG(x)\adj[\bar{E}F^\perp\bar{E}^T]=\spn\left(\{\nabla G_{12}(x), \nabla G_{22}(x)\}\right).
\] 
Then, since $\nabla G_{12}(x)=[2x_1,0]^T$ and $\nabla G_{22}(x)=[0,1]^T$, for every $x\in \R^2$, the dimension of the subspace above is 1 at $\xb$,  but it is equal to $2$ for every $x$ close enough to $\xb$ such that $x_1\neq 0$. Therefore, \sdpcr{} does not hold at $\xb$.
\jo{\hfill\Halmos \\}

\end{example}

\begin{remark}\label{rem:reduction}
Note that the construction of Lemma~\ref{lemma:convergence} does not guarantee uniqueness of $\mathcal{E}$, and $\mathcal{G}$ depends on the choices of $\bar{E}$ and also of $\mathcal{E}$. Therefore, Definition~\ref{sdp:crcq_def} also depends on $\bar{E}$ and $\mathcal{E}$. In Example~\ref{sdp:exreduction}, we considered every $\bar{E}$, and the simplest choice of $\mathcal{E}$ (constantly equal to $\bar{E}$, which is admissible for Lemma~\ref{lemma:convergence} when $r=0$), but \textit{a priori}, different $\mathcal{E}$ may lead to different variants of \sdpcr{}. Nevertheless, we should stress that they are all valid constraint qualifications and the results of this paper concerning~\sdpcr{} as in Definition~\ref{sdp:crcq_def} hold true for all choices of $\bar{E}$ and $\mathcal{E}$. Observe that the discussion of this remark, on the choice of the reduction, also applies for the NSOCP version of~\socpcr{} (Definition~\ref{socp_crcq_definition}), where a different choice of reduction function $\phi$ in \eqref{socp:reductionmap2} may lead to a different version of~\socpcr{}.
\end{remark}

Recall that the nondegeneracy condition holds at $\xb$ if, and only if, $D\mathcal{G}(\xb)\adj$ is injective. By the continuity of $D\mathcal{G}$, we have that if $\xb$ satisfies nondegeneracy, then $D\mathcal{G}(x)\adj$ remains injective for every $x$ sufficiently close to $\xb$. Then, the dimension of $D\mathcal{G}(x)\adj[F^\perp]$ remains constant for every such $x$, at every $F\faceq \S^{m-r}_+$, and it follows that nondegeneracy strictly implies \sdpcr{} as in Definition~\ref{sdp:crcq_def}.  

\begin{remark}\label{sdp_remark_eq}
Note that our approach can be trivially extended to an NSDP problem with multiple constraints. Moreover, to deal with a separate equality constraint $h(x)=0$, where $h\colon \R^n\to \R^p$, in the same spirit of Remark~\ref{socp_remark_eq}, one should consider a constraint in the form
\[
	g(x):=(G(x),h(x))
	\quad \textnormal{ and } \quad
	\K:=\S^m_+\times \{0\}^p.
\]
This yields an extension of \nlprcr{} after applying Definition~\ref{sdp:crcq_def} to the reduced form of this new problem, since $F\faceq \K$ if, and only if, $F=R\times \{0\}^p$, where $R\faceq \S^m_+$ in this case. To extend \nlpcr{} one should write the equality constraint as a pair of inequality constraints in the form $h(x)\in \R^p_+$ and $-h(x)\in\R^p_+$, giving rise to a multifold NSDP problem where $\R_+$ is seen as a copy of $\S^1_+$.
\end{remark}

For a last comparison, we should mention a constraint qualification presented in one of our previous works~\cite{seqcrcq}, which was called \emph{Sequential CRCQ} (Seq-CRCQ) therein. As a matter of fact, Seq-CRCQ differs from Definition~\ref{sdp:crcq_def} in many aspects. For instance, Example~\ref{sdp:exseq} shows that Seq-CRCQ is not implied by~\sdpcr{}. This example has already appeared in~\cite[Example 4.1]{seqcrcq}, where we show that Seq-CRCQ is not satisfied at $\xb=0$; on the other hand, we showed in Example~\ref{sdp:exseq}, that~\sdpcr{} holds at $\xb$. Thus,~\sdpcr{} is either strictly weaker than Seq-CRCQ, or completely independent of it. Despite our efforts to clarify the converse statement, we were not able to prove nor find a counterexample for it, so this is left as an open problem.

%
%
\subsection{Strong second-order optimality conditions for NSDP}\label{sec:ssocsdp}

The earliest work that provides a practical characterization of the sigma-term in NSDP is Shapiro's~\cite{shapirosdp}, using second-order directional derivatives of the least eigenvalue function, $\mu_{\min}\colon \S^m\to \R$. Shapiro proved that
\[
	\sigma(d,\xb,\bar{\lambda})=d^T \mathcal{H}(\xb,\bar{\lambda}) d,
\]
for any $d\in C(\xb)$ and $\bar{\lambda}\in \Lambda(\xb)$, where
\[
	\mathcal{H}(\xb,\bar{\lambda}):= \left[2\left\langle  
		D_{x_i}G(\xb)G(\xb)^{\dagger}D_{x_j}G(\xb) 
		 , \ 
		 \bar{\lambda}
	\right\rangle\right]_{i,j= 1,\ldots,n}
\]
and $G(\xb)^{\dag}$ denotes the Moore-Penrose pseudoinverse of $G(\xb)$. 
\if{Roughly speaking, the sigma-term can be seen as a directional derivative of a function that describes $\K$. From this point of view, Shapiro's work is motivated by the description 
\[
	\S^m_+=\{Z\in \S^m\colon \mu_{\min}(Z)\geq 0\}.
\]
Furthermore, }\fi
Shapiro also proved that if a local minimizer $\xb\in \Omega$ satisfies nondegeneracy and its associated Lagrange multiplier $\bar{\lambda}\in \Lambda(\xb)$ is such that $\rank(\bar{\lambda})+\rank(G(\xb))=m$ (strict complementarity), then $\xb$ satisfies SOC with respect to $\bar{\lambda}$.

Later, other authors provided creative ways of obtaining SOC via some local reformulation of~\eqref{nsdp} with no curvature. For instance, Lourenço et al.~\cite{louren} wrote $G(x)\succeq 0$ in the form $G(x)-Z^2=0$ with an additional variable $Z\in \S^m$, and then obtained SOC for NSDP out of SOC for NLP -- under nondegeneracy and strict complementarity. Forsgren~\cite{forsgren} rediscovered Shapiro's characterization of the sigma-term and obtained SOC (under nondegeneracy, but without strict complementarity) using a special reformulation of the problem. Jarre~\cite{Jarre} provided an elementary construction of SOC via a certain Schur complement, under nondegeneracy, strict complementarity, and assuming that the tangent cone of the linearized constraint $G(\xb)+DG(\xb)d\in \K$ coincides with $\tangent$. Fukuda et al.~\cite{fukuda2020second} used the characterization
\[
	\S^m_+=\{Z\in \S^m\mid \|\Pi_{\S^m_+}(-Z)\|^2=0\}
\]
combined with an external penalty method and the Clarke subdifferential of $\Pi_{\S^m_+}$, to achieve a weaker second-order condition, which is stated only in terms of the lineality space of $C(\xb)$. However, their results were obtained assuming only Robinson's CQ together with the so-called \textit{weak constant rank} (WCR) property, which is not a CQ on its own.

Following this line of research, the main contribution of this section consists of proving that every local minimizer $\xb$ of~\eqref{nsdp} satisfies SOC with respect to any Lagrange multiplier $\bar{\lambda}\in \Lambda(\xb)$ under the facial constant rank property. In particular, when in addition $H(\xb)$ is closed (which leads to CRCQ), then the facial constant rank property implies $\Lambda(\xb)\neq \emptyset$. \textit{A priori}, we make no special requirement on $\Lambda(\xb)$.

\begin{theorem}\label{sdp_ssoc}
Let $\xb\in \Omega$ be a local minimizer of~\eqref{nsdp} satisfying the facial constant rank property. Then, for every $\bar{\lambda}\in \Lambda(\xb)$, the pair $(\xb,\bar{\lambda})$ satisfies SOC; that is,
\begin{displaymath}
d^{T}\nabla^{2}f(\xb)d + \left\langle D^2 G(\xb)[d,d], \bar{\lambda}\right\rangle - \sigma(d,\xb,\bar{\lambda}) \geq 0
\end{displaymath}
holds for every $d\in C(\xb)=\linear\cap \{\nabla f(\xb)\}^\perp$.
\end{theorem}
\jo{\proof{Proof.}}\si{\begin{proof}}
If $\Lambda(\xb)=\emptyset$, then the result holds trivially; otherwise, let $\bar{\lambda}\in \Lambda(\xb)$ be arbitrary and fixed. Let $r:=\rank(G(\xb))$, and let $\bar{E}\in \R^{m\times m-r}$ and $\bar{P}\in \R^{m\times r}$ be matrices with orthonormal eigenvector columns associated with the zero and positive eigenvalues of $G(\xb)$, respectively. Define $\bar{U}:=[\bar{E},\bar{P}]$. Moreover, let $\mathcal{E}$ be as in Lemma~\ref{lemma:convergence} and such that $\lim_{x\to\bar{x}}E(x)=\bar{E}$, where $E(x):=\mathcal{E}(G(x))$ for every $x\in \R^n$. 

Now, let $d\in C(\xb)$ be arbitrary; so $\nabla f(\xb)^T d=0$ and $\bar{E}^T DG(\xb)d\bar{E}\succeq 0$. Following the proof of Theorem~\ref{sdp:fcrtl}, let $\bar{Q}:=[\bar{Z},\bar{W}]\in \R^{m-r\times m-r}$ be an orthogonal matrix such that $\bar{Z}^T \bar{E}^T DG(\xb)d\bar{E}\bar{Z}=0$ and $\bar{W}^T \bar{E}^T DG(\xb)d\bar{E}\bar{W}\succ 0$, and let $s$ denote the rank of $\bar{E}^T DG(\xb)d\bar{E}$. 
Moreover, let $F\faceq\S^{m-r}_+$ be defined as in~\eqref{sdp:faceproof}; that is:
\[
	F:= \left\{ \bar{Q}
		\begin{bmatrix}
			A_{11} & 0 \\ 0 & 0
		\end{bmatrix}
		\bar{Q}^T
		\ \middle| \ A_{11}\in \S^s_+	
	\right\},
\]
and note that $\bar{E}^T DG(\xb)d\bar{E}\in F$.
Similarly to the proof of Theorem~\ref{sdp:fcrtl}, since the facial constant rank property holds at $\xb$, there exists some $\varepsilon>0$ and a twice continuously differentiable curve $\xi\colon (-\varepsilon,\varepsilon)\to \R^n$ such that $\xi(0)=\xb$, $\xi'(0)=d$, and 
\[
	\mathcal{G}(\xi(t))\in \spn(F)
\]
for all $t\in (-\varepsilon,\varepsilon)$. Moreover, $\G(\xi(t))\in F$ for every $t\in[0,\varepsilon)$. Since $\xb$ is a local minimizer of~\eqref{nsdp} and $\xi(t)$ is feasible for every small $t\geq 0$, then $t=0$ is a local minimizer of the function $\phi(t):= f(\xi(t))$ subject to $t\geq 0$. Consequently, it is easy to see that
\begin{equation}\label{sdp_phi}
\phi^{\prime\prime}(0)=d^T \nabla^2 f(\xb)d+\nabla f(\xb)^T\xi^{\prime\prime}(0)\geqslant 0.
\end{equation}	

The rest of the proof consists of computing the term $\nabla f(\xb)^T\xi^{\prime\prime}(0)$. By construction, we have $\mathcal{G}(\xi(t))\in \spn(F)$ for every $t\in (-\varepsilon,\varepsilon)$, so $\bar{Z}^T \mathcal{G}(\xi(t))\bar{Z}=0$ and the reduced complementarity function
\begin{displaymath}
R(t) := \left\langle \bar{Z}^T \mathcal{G}(\xi(t))\bar{Z} , \ \bar{Z}^{T}\bar{E}^T\bar{\lambda}\bar{E}\bar{Z} \right\rangle
\end{displaymath}
has value zero, for all $t\in (-\varepsilon,\varepsilon)$. Therefore,
\begin{eqnarray*}
R^{\prime}(t) &=& \left\langle \bar{Z}^T (DE(\xi(t))\xi^{\prime}(t))^T{G}(\xi(t))E(\xi(t))\bar{Z} , \  \bar{Z}^{T}\bar{E}^T\bar{\lambda}\bar{E}\bar{Z} \right\rangle \ \\
 & & + \ \left\langle \bar{Z}^T E(\xi(t))^{T}D{G}(\xi(t))\xi^{\prime}(t)E(\xi(t))\bar{Z} , \  \bar{Z}^{T}\bar{E}^T\bar{\lambda}\bar{E}\bar{Z} \right\rangle \ \\
 & & + \ \left\langle \bar{Z}^T E(\xi(t))^{T}{G}(\xi(t))DE(\xi(t))\xi^{\prime}(t)\bar{Z} , \  \bar{Z}^{T}\bar{E}^T\bar{\lambda}\bar{E}\bar{Z} \right\rangle
\end{eqnarray*}
also has value zero for every small $t$. Differentiating once more, and taking the limit $t\to 0$, we obtain
\begin{eqnarray}\label{sdp:rterm2}
R^{\prime\prime}(0) &=& \left\langle \bar{Z}^T \bar{E}^{T}D^2{G}(\xb)[d,d] \bar{E}\bar{Z} \ + \ \bar{Z}^T \bar{E}^{T}D{G}(\xb)\xi^{\prime\prime}(0) \bar{E}\bar{Z} , \   \bar{Z}^{T}\bar{E}^T\bar{\lambda}\bar{E}\bar{Z} \right\rangle \  \nonumber\\
 & & + \ 2\left\langle \bar{Z}^T (DE(\xb)d)^T D{G}(\xb)d\bar{E}\bar{Z} \ + \ \bar{Z}^T \bar{E}^{T}D{G}(\xb)d DE(\xb)d\bar{Z} , \   \bar{Z}^{T}\bar{E}^T\bar{\lambda}\bar{E}\bar{Z} \right\rangle \  \nonumber\\
  & & + \ 2\left\langle \bar{Z}^T (DE(\xb)d)^T {G}(\xb)DE(\xb)d\bar{Z}  , \   \bar{Z}^{T}\bar{E}^T\bar{\lambda}\bar{E}\bar{Z} \right\rangle \  = \ 0. 
\end{eqnarray}
However, following Shapiro and Fan~\cite[Equation 3.8]{shapiro1995eigenvalue}, and Bonnans and Ram{\'i}rez~\cite[Equation 67]{myreportwithbonnans}, we see that
\begin{equation}\label{sdp:derivative_e}
	DE(\xb)d=D\mathcal{E}(G(\xb))DG(\xb)d=-G(\xb)^\dag DG(\xb)d \bar{E}.
\end{equation}
Substituting~\eqref{sdp:derivative_e} into~\eqref{sdp:rterm2}, the two last lines of expression~\eqref{sdp:rterm2} can be greatly simplified, which leads to the following: 
\begin{eqnarray}\label{sdp:rterm2-simplified}
R^{\prime \prime}(0) &=& \left\langle \bar{Z}^{T}\bar{E}^T\left( D^{2}{G}(\xb)[d,d] + D{G}(\xb)\xi^{\prime \prime}(0) \right)\bar{E}\bar{Z} , \   \bar{Z}^{T}\bar{E}^T\bar{\lambda}\bar{E}\bar{Z} \right\rangle \  \nonumber \\
 & & - \ \ 2\left\langle \bar{Z}^T \bar{E}^T (DG(\xb)d)^T {G}(\xb)^\dag DG(\xb)d\bar{E}\bar{Z}  , \   \bar{Z}^{T}\bar{E}^T\bar{\lambda}\bar{E}\bar{Z} \right\rangle \  = \ 0. 
\end{eqnarray}
However, by the complementarity condition we have $\bar{\lambda}\bar{P}=0$, and using the KKT condition together with $\nabla f(\xb)^T d=0$ we obtain
\[
	\begin{aligned}
		0 
		&=\langle d, -\nabla f(\xb) \rangle=\langle d, DG(\xb)\adj \bar{\lambda}\rangle=\langle DG(\xb)d, \bar{\lambda}\rangle\\
		&=\left\langle \bar{U}^T DG(\xb)d\bar{U} , \  \bar{U}^T\bar{\lambda}\bar{U}\right\rangle=\left\langle \bar{E}^T DG(\xb)d\bar{E} , \  \bar{E}^T\bar{\lambda}\bar{E}\right\rangle\\
		&=\left\langle \bar{Q}^T\bar{E}^T DG(\xb)d\bar{E}\bar{Q} , \  \bar{Q}^T\bar{E}^T\bar{\lambda}\bar{E}\bar{Q}\right\rangle\\
		&=\left\langle \bar{Z}^T\bar{E}^T DG(\xb)d\bar{E}\bar{Z} , \  \bar{Z}^T\bar{E}^T\bar{\lambda}\bar{E}\bar{Z}\right\rangle,
	\end{aligned}
\]
but since $\bar{Z}^T\bar{E}^T DG(\xb)d\bar{E}\bar{Z}\succ 0$ and $\bar{W}^T\bar{E}^T \bar{\lambda}\bar{E}\bar{W}\preceq 0$, this implies
\begin{equation*}\label{sdp:2eq1}
	\bar{Z}^T\bar{E}^T \bar{\lambda}\bar{E}\bar{Z}=0,
\end{equation*}
which in turn implies $\bar{\lambda}\bar{E}\bar{Z}=0$. With this at hand, we obtain
\begin{equation}\label{d2r}
R^{\prime \prime}(0) = \left\langle D^{2}{G}(\xb)[d,d] + D{G}(\xb)\xi^{\prime \prime}(0) - 2(DG(\xb)d)^T {G}(\xb)^\dag DG(\xb)d , \   \bar{\lambda} \right\rangle = 0,
\end{equation}
and, by the KKT conditions, this leads to
\begin{equation}\label{d2r2}
\nabla f(\xb)^T \xi^{\prime\prime}(0)=-\langle DG(\xb)\xi^{\prime\prime}(0) , \bar{\lambda} \rangle = \left\langle D^{2}{G}(\xb)[d,d] - 2(DG(\xb)d)^T {G}(\xb)^\dag DG(\xb)d , \   \bar{\lambda} \right\rangle.
\end{equation}
Substituting \eqref{d2r2} into \eqref{sdp_phi} yields
\begin{displaymath}
d^{T}\nabla^{2}f(\xb)d + \left\langle D^2 G(\xb)[d,d], \bar{\lambda}\right\rangle - d^T \mathcal{H}(\xb,\bar{\lambda})d \geq 0.
\end{displaymath}
Since $d\in C(\xb)$ was chosen arbitrarily, and $\bar{\lambda}$ is fixed from the beginning, the proof is complete. \jo{\hfill\Halmos \\}

\si{\end{proof}}

\section{Final remarks}\label{sec:conc}

The constant rank constraint qualification (CRCQ) is one of the most important regularity conditions in nonlinear programming (NLP), with several relevant applications regarding global convergence of algorithms, second-order optimality conditions, and some topics of stability theory. However, one of the main reasons why most of these interesting results still remain exclusive to NLP is that CRCQ itself seems intrinsic to NLP. Until very recently, there was no extension or analogue of it in the conic programming context. In a recent pair of papers~\cite{seqcrcq,seqcrcq-socp}, we presented an extension of CRCQ for nonlinear semidefinite and second-order cone programming using sequences and the eigenvector structure of their respective cones, which would allow us to adopt a strategy similar to the existing nonlinear programming literature. See also~\cite{weak-sparse-cq}. While this is interesting from the point of view of algorithms, it may not be an appropriate tool for other uses. Therefore, in this paper we adopted a more innovative approach: we first characterized CRCQ for NLP in a geometrical way, by means of the faces of a reduced cone, and then we showed this geometrical characterization could carry the essence of CRCQ to more general contexts. As far as we know, this is also the first time an intuitive interpretation of CRCQ was ever presented, and it is surprisingly simple: CRCQ describes the situations where every possible linear approximation of the feasible set (around a point of interest) preserves its dimension under small perturbations. As a side note, we should mention that this definition is either independent or strictly weaker than the ones presented in~\cite{seqcrcq,seqcrcq-socp}. 

As an application of our results, we obtained a strong second-order necessary optimality condition under CRCQ, in terms of any given Lagrange multiplier. This improves the classical result that is obtained under nondegeneracy, and serves as an alternative for the condition that can be obtained under Robinson's CQ, where for each direction in the critical cone, there is a Lagrange multiplier satisfying the second-order condition. We expect CRCQ to be an alternative to Robinson's CQ in other situations, especially those related with stability analysis of parametric nonlinear conic optimization programs, in view of the nonlinear programming literature -- see, for instance, references~\cite{sensitivitympec,crcq,minch}. Since CRCQ is independent of Robinson's CQ, we believe that this work allows the development of a new parallel strand in the study of stability. In a recent work, Gfrerer and Mordukhovich~\cite{tilt}, fully characterized \textit{tilt stable} local minimizers of NLP problems under the so-called \textit{bounded extreme point property}, which is implied by CRCQ (improving a previous work that assumed CRCQ and MFCQ~\cite{tiltcrmf}), and one of the possibilities of future work mentioned by them is an extension for conic programs. In a related work, still for NLP, Gfrerer and Outrata~\cite{outratagfrerer} obtained similar results to~\cite{tilt} as an application of the generalized derivative of a particular set-valued mapping, which was computed assuming the \textit{metric subregularity constraint qualification} (MSCQ) at the point of interest plus Robinson's CQ in its neighborhood. The CRCQ condition as presented in this paper may replace these assumptions in a possible extension of their results to NSOCP and NSDP. Moreover, we expect CRCQ to be useful for supporting the convergence theory of some iterative algorithms, and also to encourage the development of algorithms that rely on faces for solving nonlinear conic problems.

Regarding prospective work, the techniques employed in this paper strongly suggest that a further extension of CRCQ, for general reducible cones, is possible. In this paper we adopted a more pragmatic approach by working explicitly with NSOCP and NSDP for clarity, leaving the investigation of a more general result to future works. In fact, it would also be interesting to not rely on reducibility at all, which should be possible by taking into account the faces of the tangent cone to $\K$ at $g(\xb)$, or perturbations of it, instead of the faces of the reduced cone $\C$. Furthermore, this work may inspire extensions of weaker constant rank-type conditions from NLP (together with their applications) to the conic environment, with emphasis on the well-established \textit{constant positive linear dependence} condition~\cite{rcpld,ams05,cpld} and the \textit{constant rank of the subspace component} condition~\cite{cpg}.


\subsection*{Funding}
	This work has received financial support from
	CEPID-CeMEAI (FAPESP 2013/07375-0),
	FAPESP (grants 2018/24293-0, 2017/18308-2,
	2017/17840-2, 
	2017/12187-9, and 2020/00130-5),
	CNPq (grants
	301888/2017-5,
	303427/2018-3, and 404656/2018-8),
	PRONEX - CNPq/FAPERJ (grant E-26/010.001247/2016),
	and
	FONDECYT grant 1201982 and Basal Program CMM-AFB 170001, both from ANID (Chile).


\si{\bibliographystyle{plainurl}}
\jo{\bibliographystyle{informs2014}}
\si{
}\jo{\input{MOR-bib.bbl}}

\end{document}